\pdfoutput=1
\documentclass[11pt]{article}
\usepackage[utf8]{inputenc}
\usepackage{amsmath, amssymb, amsthm}
\usepackage{mathtools}
\usepackage{geometry}
\usepackage[expansion=false]{microtype}
\usepackage{booktabs}
\usepackage{array}
\usepackage{bm}
\usepackage{url}
\usepackage{graphicx}
\usepackage{float}
\usepackage{caption}
\usepackage{enumitem}
\usepackage{algorithm}
\usepackage{algpseudocode}
\usepackage[numbers,sort&compress]{natbib}
\usepackage{xcolor}
\definecolor{NavyBlue}{rgb}{0.0,0.0,0.5}
\usepackage[colorlinks=true,linkcolor=NavyBlue,citecolor=NavyBlue,urlcolor=NavyBlue,hypertexnames=false]{hyperref}
\hypersetup{pdftitle={Computing Multi-Slot First-Price Pacing Equilibria at Scale via Entropic Smoothing},pdfauthor={Yuchen Xie, Shenyinying Tu, Yuan Gao, Yi Zhang, Deepak Kumar}}
\usepackage[section]{placeins}
\AddToHook{cmd/subsection/before}{\FloatBarrier}

\theoremstyle{plain}
\newtheorem{theorem}{Theorem}[section]
\newtheorem{lemma}[theorem]{Lemma}
\newtheorem{proposition}[theorem]{Proposition}
\newtheorem{corollary}[theorem]{Corollary}
\theoremstyle{definition}
\newtheorem{definition}[theorem]{Definition}
\theoremstyle{remark}
\newtheorem{remark}[theorem]{Remark}

\DeclareMathOperator{\TopK}{TopK}
\newcommand{\cI}{\mathcal{I}}
\newcommand{\cJ}{\mathcal{J}}
\newcommand{\R}{\mathbb{R}}
\newcommand{\hinge}[1]{\left[#1\right]_{+}}
\newcommand{\Hb}{H_{\mathrm{bin}}}
\newcommand{\Lhard}{\mathcal{L}}
\newcommand{\Lsm}{\widetilde{\mathcal{L}}_{\mathbf{T}}}
\newcommand{\Xinf}{\mathcal{X}_{\infty}}
\newcommand{\cp}{\mathrm{cap}}

\title{Computing Multi-Slot First-Price Pacing Equilibria \\ at Scale via Entropic Smoothing\thanks{We used generative AI to help draft the text, proofs, and code; the authors reviewed all of it, and the proofs were machine-checked in Lean (Appendix~\ref{app:lean}).}}

\author{\normalsize\textbf{Yuchen Xie, Shenyinying Tu, Yuan Gao, Yi Zhang, Deepak Kumar}\\
\normalsize LinkedIn Corporation\\
\normalsize Mountain View, CA 94043\\
\normalsize\texttt{\char`\{yucxie,stu,yugao,yzhang2,dekumar\char`\}@linkedin.com}}
\date{}

\begin{document}

\maketitle

\begin{abstract}
Pacing equilibria describe how budget-constrained advertisers clear an ad market, but computing one at production scale is hard when every request fills several slots: the allocation jumps whenever two campaigns' scores cross. We consider the first-price pacing equilibrium (FPPE) of auctions with $K$ homogeneous slots and unit-demand campaigns, whose Eisenberg--Gale dual reduces, via a TopK identity, to a convex program over bids. Regularizing each request's allocation with the binary (Fermi--Dirac) entropy gives a smooth, strictly convex dual whose gradient is exactly spend minus budget and whose evaluation splits into independent one-dimensional threshold solves, one per request. The smoothed equilibrium converges to the hard one at an explicit $O(\sqrt{T})$ rate, and a duality gap certifies how far any bids are from the hard equilibrium. We implement an efficient multi-GPU solver and test it on a production ad workload. We use it to compute an equilibrium price anchor: each campaign's hindsight equilibrium cost per key result.
\end{abstract}

\section{Introduction}\label{sec:intro}

Advertising platforms allocate ad slots by auction subject to advertiser budgets, and pacing algorithms scale each campaign's bids so that its spend tracks its budget. Pacing targets a \emph{pacing equilibrium}, in which every campaign either exhausts its budget or bids at its cap. For single-slot first-price auctions the \emph{first-price pacing equilibrium} (FPPE) always exists, is unique, and solves a quasi-linear Eisenberg--Gale convex program~\citep{conitzer2022pacing,cole2017convex}. Feed and search surfaces, however, fill several slots per request, each campaign taking at most one, and \citet{conitzer2022multiplicative} posed the multi-slot extension as an open question. Concurrently and independently, \citet{afiouni2026pacing} extend the convex-program characterization to all bid-maximizing pay-your-bid mechanisms, a class that contains the multi-slot auctions studied here.

\begin{figure}[!t]
\centering
\includegraphics[width=\textwidth]{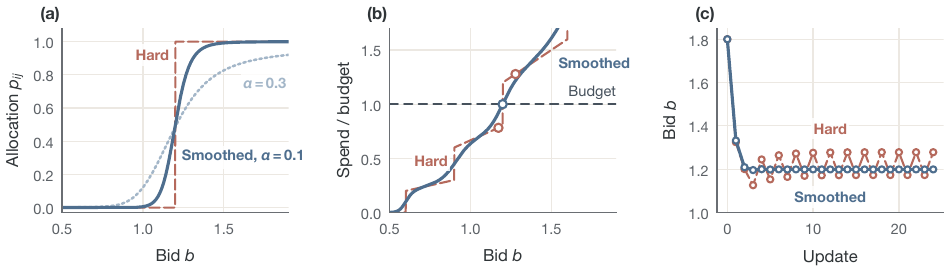}
\caption{A toy market: one campaign bids against fixed competitors in six requests. (a) Allocation in one request. (b) Spend: the hard spend jumps across the budget; the smoothed spend meets it at one bid. (c) The update of Section~\ref{sec:bid_update} without momentum: under hard top-$K$ it cycles between the circled bids in (b); under smoothing it converges.}
\label{fig:toy}
\end{figure}

This paper concerns \emph{computing} the multi-slot equilibrium at the scale of a production ad market, and the difficulty is the allocation. Selecting each request's top-$K$ winners is cheap on a GPU, but the resulting allocation is piecewise constant in the bids: its derivative vanishes almost everywhere, and it jumps wherever two campaigns' scores cross. The equilibrium solves a convex program (Section~\ref{sec:hard}), but a non-smooth one, on which first-order methods converge slowly and oscillate near score ties: a campaign at a score threshold alternately wins and loses blocks of auctions, so its spend jumps and its bid cycles instead of settling (Figure~\ref{fig:toy}). On our production workload, the update that brings every active campaign within 1\% of its budget under our smoothed allocation leaves 56\% of them further away under hard top-$K$ allocation, and a third even without momentum, with bids still cycling after 800 updates (Section~\ref{sec:exp-smoothing-ablation}).

We replace the discrete allocation with a deterministic, entropy-regularized one. The binary (Fermi--Dirac) entropy turns each request's allocation into independent Bernoulli marginals coupled by one scalar threshold, as in differentiable top-$K$ operators~\citep{amos2019lml,blondel2020fenchel,cordonnier2021difftopk}. This is not an ad hoc surrogate: the smoothed problem is the exact dual of an entropy-regularized multi-slot Eisenberg--Gale program, a smooth, strictly convex function of the log-bids whose gradient is exactly each campaign's spend minus its budget. Equilibrium computation thus remains convex minimization, now smooth and with a unique minimizer, and the smoothed equilibrium converges to the hard one as the temperature vanishes.

The same structure makes each iteration highly parallel. Given the bids, the gradient separates across requests: each request contributes the spend of its candidates and needs only one bracketed one-dimensional root solve for its threshold, and these solves run in parallel on GPUs. On this computation we build a practical solver: a log-bid update with heavy-ball momentum whose step shrinks with an annealed temperature.

\paragraph{Application: the equilibrium price anchor.} Our motivating use is a per-campaign baseline that online signals cannot provide: the cost per key result a campaign would pay \emph{in hindsight, under complete information, at pacing equilibrium}---an equilibrium cost per click for a CPC campaign, or per thousand impressions for a CPM campaign. A first-price winner pays its bid per key result, so this \textbf{equilibrium price anchor} is the campaign's equilibrium bid $b_j^\star$ (Definition~\ref{def:fppe}), well defined because equilibrium bids are unique (Theorem~\ref{thm:unique}). Realized cost per key result is confounded: it may be high because of genuine scarcity, pacing error, an infeasible budget and bid-cap combination, or audience overlap; the anchor is meant to separate market structure from controller behavior. Prospective uses include a reference for online pacing, a diagnostic for campaigns whose realized cost materially exceeds the anchor, a denoised training label for bidding and ranking models, and the evaluation of counterfactual market changes (budgets, bid caps, slot counts, or scoring rules), each one solver run. Evaluating them is beyond the scope of this paper, which concerns computing the anchor accurately and at scale, but they set the accuracy requirement: a solver that leaves even a few percent of campaigns far from their budgets leaves without a usable anchor exactly the long tail that diagnostics care about. The anchor cannot be attained online, where the pacer sees only past traffic, and we compute it on the smoothed problem, so we treat it as a smoothed proxy for the hard equilibrium (Section~\ref{sec:scope}).

\paragraph{Why first-price.}\label{sec:why_fpa} We model the market as a first-price auction because its equilibrium structure makes large-scale computation tractable. This is a modeling choice: the simulated market defines the anchor and need not reproduce the mechanism a platform deploys. Under first-price payment a campaign's spend rises with its own bid, as both its win probability and its per-win payment increase, and the pacing equilibrium solves a convex program whose smoothed gradient is exactly spend minus budget (Section~\ref{sec:fd_relaxation}). Under second-price payment the runner-up sets the price, so spend is flat in a winner's own bid over wide ranges, and pacing equilibria can be multiple and are PPAD-hard to compute~\citep{conitzer2022multiplicative,chen2024complexity}, whereas the first-price equilibrium is unique and solves a convex program~\citep{conitzer2022pacing}; Appendix~\ref{app:spa} develops the contrast.

\paragraph{Contributions.}
\begin{itemize}
    \item \textbf{Entropic smoothing with guarantees} (Sections~\ref{sec:model} and~\ref{sec:fd_relaxation}). We state the multi-slot FPPE in explicit TopK form, the homogeneous-slot case of \citet{afiouni2026pacing}. Smoothing each request's allocation by the binary entropy gives the exact dual of an entropy-regularized Eisenberg--Gale program: smooth and strictly convex, with gradient exactly spend minus budget, a sigmoid allocation tied to one bracketed threshold per request, and an explicit $O(\sqrt T)$ bound to the hard equilibrium; the same argument certifies any bids a posteriori.
    \item \textbf{A multi-GPU solver} (Section~\ref{sec:practice}). Fused Triton kernels evaluate each request in one program and shard requests across GPUs; iteration time grows with request volume and falls nearly in proportion to the number of GPUs, which does not change the results. On top runs a log-bid update with heavy-ball momentum, adaptive restart, and an annealed, temperature-coupled step; we state which parts the theory covers.
    \item \textbf{Evidence at scale} (Section~\ref{sec:empirical}). On synthetic markets solved exactly, the smoothed equilibrium approaches the exact one linearly in $T$ and the solver attains it. On a production workload the solver meets our stopping rule within 200 iterations on each of 23 days, with a certified error below 3\%; no alternative matches it: hard top-$K$ allocation, PACE, randomized budget dynamics, Adam, or plain projected gradient (Section~\ref{sec:exp-baselines}).
\end{itemize}

\section{Related Work}\label{sec:related}

\paragraph{Closest work.} Concurrently and independently, \citet{afiouni2026pacing} study pacing equilibria in abstract mechanisms. They show that the market-equilibrium interpretation of FPPE breaks down in first-price position auctions, whose ranked slots differ in quality, and identify bid-maximizing pay-your-bid mechanisms---the platform selects a feasible outcome maximizing total paced bids and charges each buyer her paced bid---as a class in which pacing equilibria exist, are unique, and solve an Eisenberg--Gale-type convex program, with Pareto-efficiency and liquid-welfare guarantees. To handle discontinuities they perturb bids randomly and multiplicatively, following \citet{borgs2007dynamics}, and prove convergence of spending-based dynamics. Our $K$-slot unit-demand auction is such a mechanism, so the hard-market results of Section~\ref{sec:hard} are the homogeneous-slot case of their theory, restated in the explicit TopK and threshold form our smoothing builds on. What we add is computational: a deterministic binary-entropy smoothing of the allocation whose dual gradient is exactly spend minus budget, with curvature bounds and a quantitative rate of convergence to the hard equilibrium (Section~\ref{sec:fd_relaxation}), and a GPU solver built on it that reaches equilibrium on a production workload (Sections~\ref{sec:practice} and~\ref{sec:empirical}).

\paragraph{Pacing and pacing equilibria.} Budget pacing in repeated ad auctions has been studied through bid dynamics~\citep{borgs2007dynamics}, fluid mean-field models~\citep{balseiro2015repeated}, online learning~\citep{balseiro2019learning}, and dual mirror descent for online allocation~\citep{balseiro2023best}; \citet{balseiro2021budget} prove that equilibria of several budget-management strategies, including multiplicative pacing, exist. \citet{gaitonde2023budget} show that gradient-based pacing attains liquid-welfare and regret guarantees without converging to an equilibrium; we compute the equilibrium itself. Pacing equilibria as a static solution concept were introduced for second-price auctions by \citet{conitzer2022multiplicative}, who show that they exist but may be multiple, and computing one is PPAD-hard~\citep{chen2024complexity}. For first-price auctions, \citet{conitzer2022pacing} prove that the pacing equilibrium is unique, coincides with the competitive equilibrium of a quasi-linear Fisher market, and solves the quasi-linear Eisenberg--Gale program of \citet{cole2017convex}; \citet{liao2023inference} build statistical inference and A/B testing on it. All these works sell one item per auction.

\paragraph{Market-equilibrium computation.} The Eisenberg--Gale program~\citep{eisenberg1959consensus} characterizes competitive equilibria of linear Fisher markets, and \citet{vazirani2007combinatorial} surveys combinatorial and convex-programming algorithms for them. \citet{birnbaum2011fisher} interpret proportional-response dynamics as mirror descent on a convex program for Fisher markets, and \citet{gao2020first} develop projected-gradient, mirror-descent, and proportional-response methods for large Fisher markets, including quasi-linear utilities. Our dual in log-bid coordinates is a multi-slot, smoothed member of this family; the smoothing makes its gradient available in closed form and its per-request computation embarrassingly parallel on a GPU.

\paragraph{Multi-slot auctions and autobidding.} Multi-slot mechanism design has a long history in sponsored search~\citep{edelman2007gsp,varian2007position}, where slots differ in click-through rate and pricing is generalized second price or VCG; \citet{ashlagi2010position} study position auctions with budget-constrained bidders from the mechanism-design side. We instead take as given a first-price rule with homogeneous slots (a campaign's value does not depend on which of the $K_i$ slots it occupies, as in our simulator) and compute the pacing equilibrium it induces; Section~\ref{sec:discussion} discusses position-dependent values. A parallel literature studies autobidders with budget and return-on-spend constraints~\citep{aggarwal2019autobidding,aggarwal2024survey} and bounds the efficiency of non-truthful and first-price auctions among them~\citep{liaw2023efficiency,deng2024efficiency}. \citet{paesleme2024complex} show that autobidding dynamics can exhibit bistability and periodic or quasi-periodic orbits instead of converging, one reason to compute the equilibrium directly rather than simulate the dynamics.

\paragraph{Differentiable top-$k$ and smoothing.} Differentiable relaxations of sorting and top-$k$ selection include NeuralSort~\citep{grover2019neuralsort}, optimal-transport operators~\citep{cuturi2019differentiable,xie2020differentiable}, regularized permutahedron projections~\citep{blondel2020fast} and their sparse top-$k$ extension~\citep{sander2023sparsetopk}, and perturbed optimizers~\citep{berthet2020learning,cordonnier2021difftopk}, with Fenchel--Young losses~\citep{blondel2020fenchel} as a common frame. The closest operator is the limited multi-label projection of \citet{amos2019lml}, which projects onto the same capped simplex with a binary-entropy regularizer and likewise obtains a sigmoid shifted by a scalar threshold found by bracketing. We use it in a different role, as the per-request smoothing inside a market-equilibrium dual: its envelope-theorem gradient is each campaign's spend and its zero-temperature limit is the hard FPPE, so we never differentiate through the root-finder as differentiable optimization layers do~\citep{amos2017optnet}. The smoothing is an instance of the frameworks of \citet{nesterov2005smooth} and \citet{beck2012smoothing}; our convergence argument avoids the general epi-convergence machinery of variational analysis~\citep{rockafellar1998variational}.

\section{Multi-Slot First-Price Pacing Equilibrium}\label{sec:model}

\subsection{Setting and notation}\label{sec:setting}

Let $\cI=\{1,\dots,N\}$ denote a batch of requests and $\cJ=\{1,\dots,M\}$ the \emph{global} set of campaigns. For each request $i$, let $\cJ_i\subseteq\cJ$ denote the eligible candidates, with $M_i:=|\cJ_i|$, typically much smaller than $M$. Request $i$ offers $K_i$ homogeneous slots, and every campaign is \emph{unit-demand}: it occupies at most one slot of any request. At most $M_i$ slots can be filled, so we assume $1\le K_i\le M_i$ without loss of generality.

For $j\in\cJ_i$, let $r_{ij}\in(0,1]$ denote the predicted rate of campaign $j$'s key result on request $i$ (the predicted click-through rate, pCTR, for CPC campaigns). Campaign $j$ has budget $B_j>0$, a bid $b_j$ per key result, and a bid cap $\cp_j$, so that $0<b_j\le\cp_j$. Its \emph{score} on request $i$ is
\begin{equation*}
    s_{ij}:=b_j\,r_{ij},
\end{equation*}
the expected amount it pays per impression if it wins a slot under first-price. We work in log-bid coordinates $x_j:=\log b_j$, write $\mathbf s_i=(s_{ij})_{j\in\cJ_i}$, and let $\Xinf:=\{\mathbf x\in\R^M: x_j\le\log\cp_j\ \forall j\}$. Per-request sums $\sum_j$ range over eligible candidates $j\in\cJ_i$; per-campaign sums $\sum_i$ range over the requests for which $j$ is eligible. Every campaign is assumed eligible for at least one request.

\paragraph{From bids to pacing multipliers.} The pacing literature writes a bid as a multiplier on value. We identify the bid cap with the campaign's value per key result---the most it is willing to pay---so that its value for request $i$ is $v_{ij}:=\cp_j\,r_{ij}$ and its \emph{pacing multiplier} is $\beta_j:=b_j/\cp_j\in(0,1]$. Then $s_{ij}=\beta_j v_{ij}$ is exactly the paced bid of \citet{conitzer2022pacing}. We state results in whichever of the equivalent coordinates $b_j$, $x_j=\log b_j$, or $\beta_j$ is most convenient.

\subsection{Background: single-slot FPPE}\label{sec:prelim}

When every request has a single slot, \citet{conitzer2022pacing} show that FPPE are exactly the optimal solutions of the quasi-linear Eisenberg--Gale program of \citet{cole2017convex}, which is program \eqref{eq:pk} below with $K_i=1$. In that program $p_{ij}$ is the (fractional) allocation of request $i$ to campaign $j$, $u_j$ is a utility proxy, and $\omega_j\ge 0$ is the budget campaign $j$ leaves unspent. The multiplier of the utility constraint is the pacing multiplier: at optimality $\beta_j=B_j/u_j\in(0,1]$, where the upper bound comes from the coefficient of $\omega_j$. Eliminating the per-request prices through $\min_{c\ge 0}\{c+\sum_j\hinge{\beta_jv_{ij}-c}\}=\max_j\beta_jv_{ij}$ leaves the reduced dual
\begin{equation}\label{eq:dql}\tag{$D_{\mathrm{QL}}$}
\min_{0<\beta_j\le 1}\ \sum_i\max_j\beta_jv_{ij}-\sum_j B_j\log\beta_j,
\end{equation}
whose unique solution gives the FPPE pacing multipliers: each campaign bids $\beta_jv_{ij}$, each request goes to its highest paced bid, and the winner pays that bid. We now derive the multi-slot analogue.

\subsection{The multi-slot Eisenberg--Gale program}\label{sec:primal}

\begin{definition}[Multi-slot unit-demand Eisenberg--Gale program]\label{def:primal}
\begin{equation}\label{eq:pk}\tag{$P_K$}
\begin{aligned}
\max_{\mathbf p,\mathbf u,\bm\omega}\quad & \sum_{j} B_j\log u_j-\sum_{j}\omega_j\\
\text{s.t.}\quad & u_j\le\sum_i v_{ij}p_{ij}+\omega_j,\quad \omega_j\ge 0 && \forall j,\\
& \textstyle\sum_{j\in\cJ_i} p_{ij}\le K_i && \forall i,\\
& 0\le p_{ij}\le 1 && \forall i,\ j\in\cJ_i.
\end{aligned}
\end{equation}
\end{definition}
Here $p_{ij}$ is the probability that campaign $j$ occupies one of request $i$'s slots; capacity limits the total mass on request $i$ to $K_i$, and $p_{ij}\le 1$ is unit demand. When $K_i=1$ capacity implies unit demand and \eqref{eq:pk} reduces to the single-slot program of \citet{cole2017convex}.

Up to the harmless bound $\omega_j\le B_j$, program \eqref{eq:pk} is the Eisenberg--Gale-type program (CP2) of \citet{afiouni2026pacing} with the feasible set taken to be the product over requests of the sets $\{\mathbf p_i\in[0,1]^{M_i}:\sum_jp_{ij}\le K_i\}$. The equilibrium notion it captures is the following.

\begin{definition}[Multi-slot FPPE]\label{def:fppe}
Bids $\mathbf b$ with $0<b_j\le\cp_j$ and an allocation $\mathbf p$ with $0\le p_{ij}\le 1$ form a \emph{multi-slot first-price pacing equilibrium} if there are thresholds $c_i\ge 0$ such that, writing $S_j:=\sum_i p_{ij}s_{ij}$ for the spend of campaign $j$,
\begin{enumerate}[label=(\roman*),leftmargin=2.2em,itemsep=1pt]
    \item \emph{Top-$K$ allocation:} $s_{ij}>c_i\Rightarrow p_{ij}=1$ and $s_{ij}<c_i\Rightarrow p_{ij}=0$;
    \item \emph{Supply:} $\sum_j p_{ij}\le K_i$, with equality whenever $c_i>0$;
    \item \emph{Budget feasibility:} $S_j\le B_j$;
    \item \emph{No unnecessary pacing:} $S_j<B_j\Rightarrow b_j=\cp_j$.
\end{enumerate}
\end{definition}
Conditions (i)--(ii) say that each request's $K_i$ slots go to its $K_i$ highest scores (ties at the threshold split fractionally) and that every winner pays its own score. For $K_i\equiv 1$ the definition coincides with the FPPE of \citet{conitzer2022pacing}.

\begin{proposition}[Feasibility, attainment, and strong duality]\label{prop:slater}
Program \eqref{eq:pk} is feasible, its optimum is attained, strong duality holds, and the Lagrangian dual optimum is attained.
\end{proposition}
The proof (Appendix~\ref{app:slater}) uses the refined Slater condition for affine constraints and a compactness argument on the superlevel sets of the objective.

\subsection{A convex dual over bids}\label{sec:hard}

\paragraph{Lagrangian dual.}\label{ssec:lagrangian}
Dualize the utility constraints with $\beta_j\ge 0$ and the capacity constraints with $c_i\ge 0$, keeping $0\le p_{ij}\le 1$ and $\omega_j\ge 0$ as domain constraints:
\begin{equation}\label{eq:lagrangian}
\mathcal{G}=\sum_j\bigl(B_j\log u_j-\omega_j\bigr)+\sum_j\beta_j\Bigl(\sum_iv_{ij}p_{ij}+\omega_j-u_j\Bigr)+\sum_i c_i\Bigl(K_i-\sum_jp_{ij}\Bigr).
\end{equation}
Maximizing over $u_j>0$ is finite only for $\beta_j>0$ and gives $B_j\log(B_j/\beta_j)-B_j$. The coefficient of $\omega_j\ge 0$ is $\beta_j-1$, so the supremum is finite only if $\beta_j\le 1$. The coefficient of $p_{ij}\in[0,1]$ is $\beta_jv_{ij}-c_i$, contributing $\hinge{\beta_jv_{ij}-c_i}$. Collecting terms, the dual is
\begin{equation}\label{eq:dk}\tag{$D_K$}
\min_{0<\beta_j\le 1,\ c_i\ge 0}\ \sum_i\Bigl[K_ic_i+\sum_j\hinge{\beta_jv_{ij}-c_i}\Bigr]-\sum_jB_j\log\beta_j+C_B,
\qquad C_B:=\sum_jB_j(\log B_j-1).
\end{equation}

\paragraph{The TopK identity.}\label{ssec:topk}
For fixed $\bm\beta$, the dual decouples across requests, and each request contributes a one-dimensional piecewise-linear problem in its threshold.

\begin{lemma}[TopK identity]\label{lem:topk}
Let $\mathbf s\in\R^n$ with $s_j\ge 0$, sorted as $s_{(1)}\ge\dots\ge s_{(n)}$, set $s_{(n+1)}:=0$, and let $1\le K\le n$. Then
\begin{equation}\label{eq:topk-identity}
\min_{c\ge 0}\Bigl\{Kc+\sum_{j=1}^n\hinge{s_j-c}\Bigr\}=\sum_{m=1}^Ks_{(m)}=:\TopK_K(\mathbf s),
\end{equation}
and the set of minimizers is exactly $[s_{(K+1)},s_{(K)}]$. If $K<n$, the same holds with the minimum taken over $c\in\R$. Moreover, $\TopK_K$ is convex, positively homogeneous, and non-decreasing in each coordinate.
\end{lemma}
The proof is in Appendix~\ref{app:topk}. For $K=1$ the identity is the $\max$ reduction that yields \eqref{eq:dql}.

\paragraph{The reduced dual.}\label{ssec:hard-dual}
Substituting Lemma~\ref{lem:topk} into \eqref{eq:dk} eliminates the thresholds. Since $\beta_jv_{ij}=s_{ij}$ and $-B_j\log\beta_j=-B_jx_j+B_j\log\cp_j$, the dual becomes a program over bids alone.

\begin{theorem}[Hard dual reduction]\label{thm:hard-dual}
Define the hard potential
\begin{equation}\label{eq:dk-hard}
\Lhard(\mathbf x):=\sum_{i}\TopK_{K_i}\bigl(\mathbf s_i(\mathbf x)\bigr)-\sum_jB_jx_j,\qquad s_{ij}(\mathbf x)=r_{ij}e^{x_j}.
\end{equation}
Then $\Lhard$ is convex, its minimum over $\Xinf$ is attained, and
\begin{equation*}
\max\eqref{eq:pk}=\min_{\mathbf x\in\Xinf}\Lhard(\mathbf x)+C_0,\qquad C_0:=\sum_jB_j\bigl(\log(B_j\cp_j)-1\bigr).
\end{equation*}
\end{theorem}
The proof is in Appendix~\ref{app:hard-dual}.

In pacing-multiplier coordinates, $\Lhard$ is the objective of the convex program (CP1) of \citet{afiouni2026pacing}, with the inner bid maximization evaluated in closed form by $\TopK$.

\paragraph{Equilibrium interpretation.} For optimal $\bm\beta^\star$, any optimal threshold $c_i^\star\in[s^\star_{i,(K_i+1)},s^\star_{i,(K_i)}]$ separates the $K_i$ winners of request $i$ from the losers: campaign $j$ wins with certainty if its score exceeds $c_i^\star$ and loses if it falls below. Each winner pays its own score (first-price pricing); the threshold is a clearing level, not a uniform price. So, unlike the single-slot case, the multi-slot FPPE is in general not a quasi-linear market equilibrium: two winners of the same request pay different prices for identical slots, and the one paying more would rather buy the cheaper slot; \citet[Theorem~1]{afiouni2026pacing} prove the analogous failure for position auctions. Nothing below relies on a market-equilibrium interpretation; the Eisenberg--Gale characterization survives without it.

\subsection{Existence and uniqueness}\label{ssec:characterization}\label{sec:unique}

\begin{theorem}[Multi-slot FPPE $\Leftrightarrow$ Eisenberg--Gale optimality]\label{thm:characterization}
$(\mathbf b,\mathbf p)$ is a multi-slot FPPE if and only if $\mathbf p$, together with some $(\mathbf u,\bm\omega)$, is optimal for \eqref{eq:pk} and $\bm\beta=\mathbf b/\bm{\cp}$, together with some $\mathbf c$, is optimal for \eqref{eq:dk}. The thresholds of Definition~\ref{def:fppe} can be taken to be the optimal dual prices. In particular, a multi-slot FPPE exists, and its bids are exactly the minimizers of $\Lhard$ over $\Xinf$ (in log coordinates).
\end{theorem}
The theorem specializes Theorem~2, Lemma~1, and Proposition~2 of \citet{afiouni2026pacing}; our proof (Appendix~\ref{app:characterization}) reads Definition~\ref{def:fppe} off the saddle-point conditions of \eqref{eq:lagrangian}: condition~(i) is maximization of $\mathcal G$ in $\mathbf p$, condition~(ii) is complementary slackness for capacity, and conditions~(iii)--(iv) follow from $u_j=B_j/\beta_j$ and $(1-\beta_j)\omega_j=0$, which give $S_j=B_j-\beta_j\omega_j$.

\begin{theorem}[Unique equilibrium bids]\label{thm:unique}
$\Lhard$ has a unique minimizer $\mathbf x^\star$ over $\Xinf$. Consequently every multi-slot FPPE has the same bids $\mathbf b^\star=e^{\mathbf x^\star}$, equivalently the same pacing multipliers $\bm\beta^\star$.
\end{theorem}
This is Corollary~1 of \citet{afiouni2026pacing} in our setting. The proof (Appendix~\ref{app:unique}) shows that the hard potential is strictly convex in bid coordinates, a fact Section~\ref{sec:convergence} reuses.

\begin{remark}
The $\TopK$ term is piecewise linear in $\mathbf b$ and contributes no curvature on each linear piece; strict convexity is carried entirely by the logarithmic Eisenberg--Gale regularizer. In log coordinates that regularizer becomes the linear term $-\sum_jB_jx_j$, which is why uniqueness is proved in bid coordinates.
\end{remark}

\begin{corollary}[Essential uniqueness]\label{cor:essential}
In every multi-slot FPPE: (a) the bids are $\mathbf b^\star$; (b) every paced campaign ($b^\star_j<\cp_j$) spends exactly $B_j$; (c) the revenue of request $i$, $\sum_jp_{ij}s^\star_{ij}$, equals $\TopK_{K_i}(\mathbf s_i^\star)$; and (d) at every request with $s^\star_{i,(K_i)}>s^\star_{i,(K_i+1)}$ the allocation is unique: the $K_i$ highest-score campaigns win with probability one. The allocation can differ across equilibria only among campaigns tied at the boundary score of a request.
\end{corollary}
\citet{afiouni2026pacing} observe that total revenue is the same for every bid-maximizing allocation at the equilibrium multipliers; part (c) refines this to every request.
The proof is in Appendix~\ref{app:essential}.

\section{Entropic Smoothing and Its Guarantees}\label{sec:fd_relaxation}

The hard potential $\Lhard$ is convex but not smooth. $\TopK_{K_i}$ has kinks wherever the $K_i$-th and $(K_i{+}1)$-th scores cross, and its gradient is the $0/1$ winner indicator, so the gradient of $\Lhard$ is constant between crossings and jumps across them; replaying the discrete auction inherits exactly this behavior. We therefore smooth the per-request allocation so that the dual stays a convex program over bids, now smooth and strictly convex (Proposition~\ref{prop:convex_smooth}). Table~\ref{tab:comparison} summarizes how each ingredient generalizes its single-slot counterpart.

\begin{table}[htbp]
\centering
\small
\caption{Single-slot versus multi-slot first-price pacing. With $K_i=1$ the hard quantities coincide; the smoothed operators differ but have the same hard limit.}
\label{tab:comparison}
\renewcommand{\arraystretch}{1.2}
\begin{tabular}{@{}p{0.24\textwidth}p{0.28\textwidth}p{0.42\textwidth}@{}}
\toprule
 & \textbf{Single slot} ($K_i=1$) & \textbf{Multi-slot} (this paper)\\
\midrule
Per-request dual operator & $\max_js_{ij}$ & $\TopK_{K_i}(\mathbf s_i)$\\
Winners per request & 1 & $K_i$, at most one slot each\\
Unit demand & implied by supply & explicit, $p_{ij}\le 1$\\
Dual non-smoothness & $s_{i,(1)}=s_{i,(2)}$ & $s_{i,(K_i)}=s_{i,(K_i+1)}$\\
Entropic regularizer & Shannon, on the simplex & binary (Fermi--Dirac), on the capped simplex\\
Smoothed value & $T\log\sum_je^{s_{ij}/T}$ & $\min_c\{K_ic+T\sum_j\log(1+e^{(s_{ij}-c)/T})\}$\\
Soft allocation & softmax $\pi_{ij}$ & root-tied sigmoid $p_{ij}$\\
Smoothed gradient & $\sum_i\pi_{ij}s_{ij}-B_j$ & $\sum_ip_{ij}s_{ij}-B_j$\\
Inner computation & closed form & one monotone 1-D root per request\\
\bottomrule
\end{tabular}
\end{table}

\paragraph{Why not softmax.} In a single-slot auction the classical continuous relaxation replaces the winner indicator by a softmax $\pi_{ij}\propto\exp(s_{ij}/T)$ with temperature $T>0$: over the simplex it maximizes expected payment plus $T$ times the Shannon entropy, and its value is the log-sum-exp smoothing of $\max_js_{ij}$. Softmax does not extend to multi-slot auctions: its coupled denominator ties every candidate's allocation to every other's, so five identical top bids each get mass $1/5$, never $1$, and the rescaled $K_i\cdot\pi_{ij}$ can put more than one unit on a single campaign, violating unit demand. It structurally cannot allocate full mass to multiple winners.

\subsection{Fermi--Dirac regularization of the allocation}

The feasible allocation space for a single request $i$ is the capped $K_i$-simplex
\begin{equation*}
    \mathcal{C}_{K_i}=\Big\{\mathbf p_i\in\R^{M_i}\ \Big|\ 0\le p_{ij}\le 1\ \forall j\in\cJ_i,\ \sum_{j\in\cJ_i}p_{ij}=K_i\Big\}.
\end{equation*}
Each request $i$ carries its own temperature $T_i>0$; we write $\mathbf T=(T_1,\ldots,T_N)$ for the whole collection and elide the subscript where only one request is in view. We maximize the expected payment regularized by the Fermi--Dirac (binary) entropy, which penalizes each marginal separately and keeps it in $[0,1]$---exactly the unit-demand structure:
\begin{equation}\label{eq:primal}
    V_i(\mathbf s_i;T_i):=\max_{\mathbf p_i\in\mathcal C_{K_i}}\ \sum_{j}p_{ij}\,s_{ij}+T_i\sum_{j}\Hb(p_{ij}),\qquad \Hb(p)=-p\log p-(1-p)\log(1-p).
\end{equation}

\begin{lemma}[Threshold representation]\label{lem:fd-dual}
Let $K_i<M_i$ and $T_i>0$. Then
\begin{equation}\label{eq:Vi}
    V_i(\mathbf s_i;T_i)=\min_{c\in\R}\Bigl\{K_i\,c+T_i\sum_{j}\log\Bigl(1+\exp\frac{s_{ij}-c}{T_i}\Bigr)\Bigr\}.
\end{equation}
The minimum is attained at a unique $c_i$, and the unique maximizer of \eqref{eq:primal} is
\begin{equation}\label{eq:sigmoid}
    p_{ij}=\sigma\!\left(\frac{s_{ij}-c_i}{T_i}\right),\qquad \sigma(z)=\frac{1}{1+e^{-z}}.
\end{equation}
If $K_i=M_i$, then $\mathcal C_{K_i}=\{\mathbf 1\}$, $V_i(\mathbf s_i;T_i)=\sum_js_{ij}$, and $p_{ij}=1$.
\end{lemma}
The proof is in Appendix~\ref{app:fd}. It rests on the per-coordinate identity $\max_{p\in[0,1]}\{pz+T\Hb(p)\}=T\log(1+e^{z/T})$, attained at $p=\sigma(z/T)$; the right-hand side, the antiderivative of the logistic function, is sometimes called the \emph{softplus}; we write it out explicitly throughout. The two descriptions of $V_i$ are Fenchel dual to each other: \eqref{eq:primal} regularizes the allocation, while \eqref{eq:Vi} is the threshold form of $\TopK$ from Lemma~\ref{lem:topk} with each hinge $\hinge{z}$ replaced by its softplus $T_i\log(1+e^{z/T_i})$.

Unlike softmax, \eqref{eq:sigmoid} \emph{decouples} the candidate denominators: each $p_{ij}$ is a Bernoulli over its own score relative to a single ``ghost bid'' $c_i$, and only the scalar capacity constraint ties the request together. As $T_i\to 0^+$, each sigmoid steepens into a step; whenever the $K_i$-th and $(K_i{+}1)$-th scores are strictly separated, the root-finder slides $c_i$ into the gap between them and recovers exact hard top-$K$ (Lemma~\ref{lem:mismass} below makes this quantitative). When scores tie at the threshold, the hard top-$K$ solution is itself non-unique: the strictly higher-scored candidates plus any $K_i-(\text{above-tie count})$ of the tied boundary candidates give an equal-value hard solution. The relaxation picks a canonical soft split among the tied candidates, so it resolves the hard-tie ambiguity rather than introducing it. (At our final temperature, 18.2\% of audited requests place more than $0.01$ of their mass outside the top $K_i$, all of them with a boundary gap of at most 11.4 temperatures; see Appendix~\ref{sec:kernel_correctness}.)

\begin{lemma}[Smoothing gap]\label{lem:gap}
For every request $i$ and score vector $\mathbf s_i$,
\begin{equation}\label{eq:smoothing_gap}
\begin{split}
    \TopK_{K_i}(\mathbf s_i)\ \le\ V_i(\mathbf s_i;T_i)\ &\le\ \TopK_{K_i}(\mathbf s_i)+T_i\,M_i\,\Hb\!\left(\tfrac{K_i}{M_i}\right)\\
    &\le\ \TopK_{K_i}(\mathbf s_i)+T_i\min\Bigl\{M_i\log 2,\ K_i\bigl(1+\log\tfrac{M_i}{K_i}\bigr)\Bigr\}.
\end{split}
\end{equation}
\end{lemma}
The lower bound holds because the hard top-$K_i$ indicator is feasible in \eqref{eq:primal} and has zero entropy; the upper bound follows from Jensen's inequality applied to the concave $\Hb$ (Appendix~\ref{app:gap}).

\subsection{A regularized market with an exact dual}

Working in log-bid space, we form the smoothed Eisenberg--Gale dual across all requests, penalized by target budgets $B_j$:
\begin{equation}\label{eq:dual}
    \Lsm(\mathbf x):=\sum_{i=1}^NV_i\bigl(\mathbf s_i(\mathbf x);\,T_i\bigr)-\sum_{j=1}^MB_j\,x_j,\qquad s_{ij}(\mathbf x)=r_{ij}e^{x_j},
\end{equation}
where the subscript records the per-request temperatures. $\Lsm$ is not merely a smooth surrogate for $\Lhard$: it is the exact dual of an entropy-regularized Eisenberg--Gale program whose optimal allocation is the Fermi--Dirac soft top-$K$.

\begin{proposition}[Regularized Eisenberg--Gale program]\label{prop:reg-market}
Let $(P_K^{\mathbf T})$ denote program \eqref{eq:pk} with each capacity constraint imposed with equality and the regularizer $\sum_iT_i\sum_{j}\Hb(p_{ij})$ added to the objective. Then
\begin{equation*}
    \max\,(P_K^{\mathbf T})=\min_{\mathbf x\in\Xinf}\Lsm(\mathbf x)+C_0,
\end{equation*}
both optima are attained, and for any minimizer $\mathbf x$ of $\Lsm$ the optimal allocation is the sigmoid \eqref{eq:sigmoid} evaluated at $\mathbf s_i(\mathbf x)$. With $S_j:=\sum_ip_{ij}s_{ij}$, the pair $(\mathbf b,\mathbf p)=(e^{\mathbf x},\mathbf p)$ satisfies budget feasibility and no unnecessary pacing (Definition~\ref{def:fppe}(iii)--(iv)). We call it the \emph{smoothed FPPE} at temperature $\mathbf T$.
\end{proposition}
The proof is in Appendix~\ref{app:fd}. The anchor computed from the smoothed problem is therefore itself an equilibrium quantity: the equilibrium bid of the mechanism in which each request allocates its slots by the soft top-$K$ rule \eqref{eq:sigmoid} and every winner pays its own score.

\subsection{One threshold per request}\label{ssec:sigmoid}

The threshold $c_i$ in \eqref{eq:sigmoid} is determined by the capacity constraint.

\begin{theorem}[Root-tied sigmoid]\label{thm:sigmoid}
Let $K_i<M_i$ and $T_i>0$. The function $F_i(c):=\sum_j\sigma\bigl((s_{ij}-c)/T_i\bigr)$ is continuous and strictly decreasing, with $F_i'(c)=-\frac{1}{T_i}\sum_jp_{ij}(1-p_{ij})<0$, and maps $\R$ onto $(0,M_i)$. Hence the root equation
\begin{equation}\label{eq:sigmoid-root}
    F_i(c_i)=K_i
\end{equation}
has a unique solution $c_i^\star$, which is a $C^\infty$ function of $(\mathbf s_i,T_i)$ and is the minimizer in \eqref{eq:Vi}. With $\rho_i:=T_i\log\frac{M_i-K_i}{K_i}$,
\begin{equation}\label{eq:bracket}
    \min_js_{ij}+\rho_i\ \le\ c_i^\star\ \le\ \max_js_{ij}+\rho_i.
\end{equation}
\end{theorem}
The proof is in Appendix~\ref{app:sigmoid}.

\begin{corollary}[Uniformly valid symmetric brackets]\label{cor:bracket}
The bracket $[\min_js_{ij}-\kappa T_i,\ \max_js_{ij}+\kappa T_i]$ contains the root for every score vector if and only if $\kappa\ge\bigl|\log\frac{M_i-K_i}{K_i}\bigr|$.
\end{corollary}
Sufficiency is \eqref{eq:bracket}; necessity follows by taking all scores equal, in which case the root is exactly $s+\rho_i$ (Appendix~\ref{app:sigmoid}).

\begin{lemma}[Mis-assigned mass]\label{lem:mismass}
Let $K_i<M_i$ and suppose the boundary scores are separated by a gap $g_i:=s_{i,(K_i)}-s_{i,(K_i+1)}>0$. Then the soft allocation places mass $U_i:=\sum_{j\notin\text{top-}K_i}p_{ij}$ outside the top $K_i$ candidates and leaves mass $D_i:=\sum_{j\in\text{top-}K_i}(1-p_{ij})$ missing from them, with
\begin{equation*}
    U_i=D_i\le\max\{K_i,\,M_i-K_i\}\,e^{-g_i/(2T_i)}.
\end{equation*}
\end{lemma}
The proof is in Appendix~\ref{app:mismass}.

\subsection{The gradient is spend minus budget}

\begin{theorem}[Budget-gap gradient identity]\label{thm:gradient}
For fixed $\mathbf T$, $\Lsm$ is differentiable. Using $\frac{d}{dz}\log(1+e^z)=\sigma(z)$ and $\frac{\partial s_{ij}}{\partial x_j}=s_{ij}$, the direct dependence gives
\begin{equation}\label{eq:grad}
    \frac{\partial V_i}{\partial x_j}=\sigma\!\left(\frac{s_{ij}-c_i}{T_i}\right)\cdot\frac{\partial s_{ij}}{\partial x_j}=p_{ij}\,s_{ij},
\end{equation}
the expected spend contribution of campaign $j$ on request $i$, and summing over requests,
\begin{equation}\label{eq:SminusB}
    \boxed{\;\nabla_{x_j}\Lsm(\mathbf x)=S_j(\mathbf x;\mathbf T)-B_j,\qquad S_j(\mathbf x;\mathbf T):=\sum_ip_{ij}\,s_{ij}.\;}
\end{equation}
\end{theorem}
The proof (Appendix~\ref{app:gradient}) applies the envelope theorem to the threshold representation \eqref{eq:Vi}: since $c_i$ is an inner minimizer, its dependence on the scores drops out.

The rest of the paper builds on this identity: equilibrium-finding needs no finite differences, no autograd through the root-finder, and no stochastic gradient estimators. A descent step on $\Lsm$ lowers the bids of campaigns that overspend and raises the bids of those that underspend: it is the classical pacing fixed-point iteration, now with an exact potential behind it.

\subsection{Convexity, smoothness, and curvature}

Write $a_{ij}:=p_{ij}(1-p_{ij})$ for the Bernoulli variance of candidate $j$ on request $i$ and $w_{ik}:=a_{ik}/\sum_\ell a_{i\ell}$ for its normalized weight (for requests with $K_i=M_i$, all $a_{ij}=0$ and we set $w_{ik}:=0$). Differentiating the root condition $\sum_jp_{ij}=K_i$ gives $\partial c_i/\partial s_{ik}=w_{ik}\in[0,1]$ and hence
\begin{equation}\label{eq:dpds}
    \frac{\partial p_{ij}}{\partial s_{ik}}=\frac{1}{T_i}\,a_{ij}\,\bigl[\delta_{jk}-w_{ik}\bigr].
\end{equation}

\begin{proposition}\label{prop:convex_smooth}
Fix a temperature vector $\mathbf T$ with $T_i>0$ for all $i$, and assume $r_{ij}>0$, $1\le K_i\le M_i$, and $B_j>0$. Then the dual $\Lsm:\R^M\to\R$ defined in \eqref{eq:dual} is:
\begin{enumerate}
    \item[(i)] convex in $\mathbf x$; more precisely, for every $\mathbf x$,
    \begin{equation}\label{eq:hessian-sandwich}
        \mathrm{diag}\bigl(S_j\bigr)\ \preceq\ \nabla^2\Lsm(\mathbf x)\ \preceq\ \mathrm{diag}\Bigl(S_j+\tfrac14\textstyle\sum_is_{ij}^2/T_i\Bigr),
    \end{equation}
    so $\Lsm$ is strictly convex;
    \item[(ii)] $C^\infty$-smooth in $\mathbf x$;
    \item[(iii)] gradient-Lipschitz on every compact convex set $\mathcal Q$, with Lipschitz constant
    \begin{equation*}
        L_\mathbf T\ \le\ \max_{\mathbf x\in\mathcal Q}\max_j\Bigl(S_j+\tfrac14\textstyle\sum_is_{ij}^2/T_i\Bigr)\ =\ O\!\left(\max_i 1/T_i\right).
    \end{equation*}
\end{enumerate}
\end{proposition}

The proof (Appendix~\ref{app:convex}) is per-request, so per-request temperatures do not weaken the result; temperatures that depend on $\mathbf x$ would, and \S\ref{sec:scope} addresses them. Convexity follows because $V_i$ is a supremum of affine functions of the scores and is non-decreasing in each, composed with the convex map $x_j\mapsto r_{ij}e^{x_j}$; the Hessian bounds follow from \eqref{eq:dpds} and the Cauchy--Schwarz inequality.

\emph{Remark on the Lipschitz rate.} The $O(1/T)$ behavior is tight: as $T\to 0^+$, the sigmoid $\sigma((s-c)/T)$ approaches a step function, and finite differences of $S_j$ across the threshold become unbounded. This quantifies the intuition of sharp gradients near the decision boundary and motivates the $\alpha_{\text{step}}\propto T$ coupling of Section~\ref{sec:bid_update}.

The Hessian formula also makes the first-price intuition of Section~\ref{sec:why_fpa} precise.

\begin{proposition}[First-price spend elasticities]\label{prop:elasticity}
For fixed $\mathbf T$ and every $\mathbf x$, the own-bid and cross-bid spend responses are
\begin{equation*}
    \frac{\partial S_j}{\partial x_j}=S_j+\sum_i\frac{a_{ij}s_{ij}^2}{T_i}\bigl(1-w_{ij}\bigr)\ \ge\ S_j,
    \qquad
    \frac{\partial S_j}{\partial x_k}=-\sum_i\frac{s_{ij}s_{ik}}{T_i}\,a_{ij}w_{ik}\ \le\ 0\quad(k\ne j).
\end{equation*}
Hence the own-bid spend elasticity satisfies $\eta_j:=\partial\log S_j/\partial\log b_j\ge 1$, and raising a competitor's bid never raises campaign $j$'s spend.
\end{proposition}
The proof differentiates $S_j=\sum_ip_{ij}s_{ij}$ using \eqref{eq:dpds} (Appendix~\ref{app:convex}). The first term, $S_j$, is the price effect: under first-price payment a 1\% higher bid raises the payment on every won impression by 1\%; the second is the win-probability effect. Under second-price payment the price effect is absent, which is the plateau of Appendix~\ref{app:spa}(b). Integrating the first inequality, $S_j$ grows at least proportionally with $b_j$ when the other bids are held fixed: this is the own-bid responsiveness condition (Assumption~5) of \citet{afiouni2026pacing} with exponent $\mu_j=1$, and the sign of the cross effect is their opponent-inverse-monotonicity (payment-monotonicity) condition (Assumption~1).

\subsection{Uniqueness at every temperature}

\begin{proposition}[Unique smoothed equilibrium]\label{prop:smooth-unique}
For every $\mathbf T>0$, $\Lsm$ has a unique minimizer $\mathbf x^\star(\mathbf T)$ over $\Xinf$ (and over $\mathcal X$).
\end{proposition}
The proof is in Appendix~\ref{app:smooth-unique}.

\subsection{Convergence to the hard equilibrium}\label{sec:convergence}

\begin{theorem}[Convergence under annealing]\label{thm:convergence}
Let $\mathbf b^\star$ be the unique multi-slot FPPE bids (Theorem~\ref{thm:unique}) and $\mathbf b^\star(\mathbf T)=e^{\mathbf x^\star(\mathbf T)}$ the unique smoothed equilibrium bids (Proposition~\ref{prop:smooth-unique}). Then
\begin{equation}\label{eq:rate}
    \sum_jB_j\left(\frac{b^\star_j(\mathbf T)-b^\star_j}{\cp_j}\right)^2=\sum_jB_j\bigl(\beta^\star_j(\mathbf T)-\beta^\star_j\bigr)^2\ \le\ 4\,\Gamma_\mathbf T,
\end{equation}
where
\begin{equation*}
    \Gamma_\mathbf T:=\sum_iT_iM_i\Hb\!\left(\tfrac{K_i}{M_i}\right)\ \le\ \sum_iT_iK_i\Bigl(1+\log\tfrac{M_i}{K_i}\Bigr).
\end{equation*}
In particular $\mathbf b^\star(\mathbf T)\to\mathbf b^\star$ as $\max_iT_i\to 0^+$, with error $O\bigl(\sqrt{\max_iT_i}\bigr)$.
\end{theorem}
\begin{proof}[Proof sketch]
Three facts combine; the full proof is in Appendix~\ref{app:convergence}. (1)~\emph{One-sided approximation:} by Lemma~\ref{lem:gap}, $0\le\Lsm-\Lhard\le\Gamma_\mathbf T$ uniformly. (2)~\emph{Near-optimality:} hence $\Lhard(\mathbf x^\star(\mathbf T))\le\Lsm(\mathbf x^\star(\mathbf T))\le\Lsm(\mathbf x^\star)\le\Lhard(\mathbf x^\star)+\Gamma_\mathbf T$. (3)~\emph{Strong convexity in bid coordinates:} on the box $0<b_j\le\cp_j$, the regularizer $-\sum_jB_j\log b_j$ has curvature at least $B_j/\cp_j^2$, so $\Lhard(\log\mathbf b)$ is $1$-strongly convex in the norm $\|\mathbf y\|_W^2=\sum_j(B_j/\cp_j^2)y_j^2$, and a midpoint argument turns a $\Gamma_\mathbf T$-minimizer into a point within $W$-distance $2\sqrt{\Gamma_\mathbf T}$ of $\mathbf b^\star$.
\end{proof}

The bound is worst-case, and it concerns the \emph{exact} relaxed minimizer $\mathbf x^\star(\mathbf T)$, not the solver's iterate (Section~\ref{sec:scope}). Appendix~\ref{app:convergence} discusses its scale and shows that, at every request whose FPPE boundary scores are separated, the soft allocation also converges to the hard one.

The same strong convexity also bounds the distance of \emph{any} bids from the hard equilibrium, by a duality gap that can be evaluated at the solver's output.

\begin{corollary}[A posteriori certificate]\label{cor:certificate}
Let $\mathbf x\in\Xinf$, write $\beta_j=e^{x_j}/\cp_j$, and let $\mathbf p$ be any allocation feasible for \eqref{eq:pk}, with values $V_j:=\sum_iv_{ij}p_{ij}$. Define
\begin{equation*}
    \Delta(\mathbf x,\mathbf p):=\Lhard(\mathbf x)+C_0-\sum_j\Bigl(B_j\log\max\{V_j,B_j\}-\hinge{B_j-V_j}\Bigr).
\end{equation*}
Then $\Delta(\mathbf x,\mathbf p)\ge\Lhard(\mathbf x)-\min_{\Xinf}\Lhard\ge0$ and
\begin{equation*}
    \sum_jB_j\bigl(\beta_j-\beta^\star_j\bigr)^2\ \le\ 4\,\Delta(\mathbf x,\mathbf p).
\end{equation*}
\end{corollary}
\begin{proof}
For fixed $\mathbf p$, the best $(\mathbf u,\bm\omega)$ in \eqref{eq:pk} is $\omega_j=\hinge{B_j-V_j}$ and $u_j=V_j+\omega_j$, with objective value $\sum_j(B_j\log\max\{V_j,B_j\}-\hinge{B_j-V_j})$. This is at most $\max\eqref{eq:pk}=\min_{\Xinf}\Lhard+C_0$ (Theorem~\ref{thm:hard-dual}), which gives the first inequality. Step~3 of the proof of Theorem~\ref{thm:convergence} (Appendix~\ref{app:convergence}) bounds $\sum_jB_j(\beta_j-\beta^\star_j)^2$ by four times $\Lhard(\mathbf x)-\min_{\Xinf}\Lhard$.
\end{proof}
Under first-price payment $V_j=S_j(\mathbf p)/\beta_j$, where $S_j(\mathbf p)=\sum_ip_{ij}s_{ij}$ is the spend of campaign $j$ under $\mathbf p$, and $\Lhard(\mathbf x)$ needs one hard top-$K$ pass at $\mathbf x$. The certificate therefore costs two auction passes at the final bids, one hard and one smoothed, and it applies to the solver's actual output whatever the departures (S1)--(S3) of Section~\ref{sec:scope}. Section~\ref{sec:convergence_table} evaluates it on the production workload.

\section{A Practical Solver}\label{sec:practice}

We first state a reference algorithm whose guarantees follow from Section~\ref{sec:fd_relaxation}, then the practical solver used in our experiments, and finally how it departs from the theory.

\subsection{A reference algorithm}\label{sec:reference}

With $\nabla_{x_j}\Lsm=S_j-B_j$ from \eqref{eq:SminusB}, equilibrium-finding at fixed $\mathbf T$ is smooth, strictly convex minimization on a box. Our practical solver is not a descent method on any single $\Lsm$ but a log-space update with momentum and an outer temperature anneal, motivated by, but not formally derived from, this convex structure; the reference method is the one the theorems back.

Algorithm~\ref{alg:solver} is annealed projected gradient descent on $\Lsm$. At fixed $\mathbf T$, with a constant step below $1/L_\mathbf T$ on the initial sublevel set it converges to $\mathbf x^\star(\mathbf T)$, and linearly because $\nabla^2\Lsm\succeq\mathrm{diag}(\mathbf S)$ is positive definite on that compact set (Proposition~\ref{prop:convex_smooth}; see~\citep{bertsekas1999nonlinear,nocedal2006numerical}). Every step is two GPU kernels: a per-request root solve and a per-campaign scatter-add. In pacing-multiplier coordinates the update is $\beta_j\leftarrow\min\{\beta_je^{-\lambda_t(S_j-B_j)},1\}$, the multiplicative-weights pacing dynamics that \citet{afiouni2026pacing} analyze in time average for perturbed mechanisms; by Theorem~\ref{thm:gradient}, for the Fermi--Dirac-smoothed first-price mechanism it is projected gradient descent on an explicit smooth convex potential, so at fixed $\mathbf T$ the iterates themselves converge.

\begin{algorithm}[htbp]
\caption{Annealed projected gradient on the smoothed dual (reference solver)}
\label{alg:solver}
\begin{algorithmic}[1]
\Require pCTRs $\{r_{ij}\}$, budgets $\{B_j\}$, caps $\{\cp_j\}$, slot counts $\{K_i\}$, initial temperatures $\{T_i^{(0)}\}$, anneal factor $\theta\in(0,1]$, step sizes $\{\lambda_t\}$, tolerance $\epsilon_{\mathrm{tol}}$
\State $x_j\gets\log\cp_j$ for all $j$ \Comment{start unpaced: $\beta_j=1$}
\Repeat
  \For{each request $i$ \textbf{in parallel}}
    \State $s_{ij}\gets r_{ij}e^{x_j}$ for $j\in\cJ_i$
    \If{$K_i=M_i$}
      \State $p_{ij}\gets 1$ for $j\in\cJ_i$ \Comment{exact (Lemma~\ref{lem:fd-dual})}
    \Else
      \State $\rho_i\gets T_i\log\bigl((M_i-K_i)/K_i\bigr)$
      \State solve $\sum_j\sigma\bigl((s_{ij}-c_i)/T_i\bigr)=K_i$ for $c_i\in[\min_js_{ij}+\rho_i,\ \max_js_{ij}+\rho_i]$ \Comment{Theorem~\ref{thm:sigmoid}}
      \State $p_{ij}\gets\sigma\bigl((s_{ij}-c_i)/T_i\bigr)$ for $j\in\cJ_i$
    \EndIf
  \EndFor
  \State $S_j\gets\sum_ip_{ij}s_{ij}$ for all $j$ \Comment{scatter-add; $\nabla_{x_j}\Lsm=S_j-B_j$}
  \State $x_j\gets\min\bigl\{x_j-\lambda_t(S_j-B_j),\ \log\cp_j\bigr\}$ for all $j$ \Comment{projected gradient step}
  \State $T_i\gets\theta\,T_i$ for all $i$ \Comment{anneal}
\Until{the KKT residual of \S\ref{sec:kkt} is at most $\epsilon_{\mathrm{tol}}$}
\Ensure bids $b_j=e^{x_j}$, soft allocation $\{p_{ij}\}$
\end{algorithmic}
\end{algorithm}

Our practical update (Section~\ref{sec:bid_update}) differs in two ways: it steps along a clipped log budget-to-spend ratio rather than along $S_j-B_j$, which divides the gradient by each campaign's budget to first order, and it adds momentum. With budgets that span orders of magnitude, no single scalar step serves Algorithm~\ref{alg:solver} well (Section~\ref{sec:exp-baselines}).

\subsection{Optimality on a bounded box}\label{sec:kkt}

The solver restricts $\mathbf x$ to the compact box $\mathcal X=\{\mathbf x:\log\varepsilon\le x_j\le\log\cp_j\}$, on which the continuous $\Lsm$ (Proposition~\ref{prop:convex_smooth}(ii)) attains its minimum. The first-order condition on the box is the projected-gradient condition, \emph{not} $S_j=B_j$ for every $j$: $S_j=B_j$ for an interior coordinate ($\log\varepsilon<x_j^\star<\log\cp_j$), $S_j\le B_j$ for a capped campaign ($x_j^\star=\log\cp_j$; it would bid more but cannot), and $S_j\ge B_j$ for a floored one ($x_j^\star=\log\varepsilon$). A convergence criterion must therefore partition campaigns by active set: a capped campaign with $S_j<B_j$ is converged. The floor is an implementation device, not part of the economic model, and it is inactive once $\varepsilon$ is small enough.

\begin{lemma}[Inactive floor]\label{lem:floor}
Let $R_j:=\sum_ir_{ij}$ and $\underline b_j:=\min\{\cp_j,\,B_j/R_j\}$. Every minimizer of $\Lhard$ or of $\Lsm$ over $\Xinf$ satisfies $b_j\ge\underline b_j$ for all $j$. Hence if $\varepsilon\le\min_j\underline b_j$, minimizing over $\mathcal X$ or over $\Xinf$ yields the same unique minimizer.
\end{lemma}
The proof is in Appendix~\ref{app:floor}: below the cap a minimizer must spend its budget, and spend is at most $b_jR_j$. An analogous bound on their pacing multipliers, $\ge\min\{1,B_i/\bar V\}$ with $\bar V$ a bound on valuations, appears in the analysis of PACE-style dynamics of \citet{afiouni2026pacing}.

\subsection{Computing the thresholds on a GPU}\label{sec:solver}

Computing $p_{ij}$ requires the root $c_i$ of \eqref{eq:sigmoid-root} for every request. Requests are independent given the bids, so all roots are solved in parallel by a vectorized kernel with no per-request loop. Each root is bracketed by $[s_{\min}-\kappa T_i,\ s_{\max}+\kappa T_i]$, where $s_{\min},s_{\max}$ are the request's smallest and largest scores; by Corollary~\ref{cor:bracket} such a bracket is valid for every score vector exactly when $\kappa\ge|\log((M_i-K_i)/K_i)|$, which $\kappa=50$ satisfies for any realistic candidate pool. A fixed number of bisection steps followed by Newton steps safeguarded to stay inside the bracket converges quadratically near the root. Appendix~\ref{app:numerics} gives the iteration counts and the safeguards needed in \texttt{float32}.

\paragraph{Fused multi-GPU implementation.} The main configuration evaluates each auction pass with fused Triton kernels~\citep{tillet2019triton}. One program handles one request: it loads the key-result rates of the request's candidates and the current bids, forms the scores, computes the temperature from their spread, solves for the threshold, forms the soft allocation and rescales it to capacity, and adds each candidate's expected payment to its campaign's spend by atomic accumulation in double precision; no per-candidate intermediate array is written to device memory. Requests are routed by structure: when every candidate wins ($K_i=M_i$) the allocation is exact and needs no root (Lemma~\ref{lem:fd-dual}); requests with two candidates or a single slot use closed-form or specialized roots; all others use the bracketed bisection--Newton solver above, register-resident for small candidate pools and tiled for large ones. A root that misses its tolerance is recomputed in double precision, and a pass with any unresolved request is rejected rather than used for a bid update. Requests are sharded across GPUs, which evaluate their shards concurrently, and the per-campaign spends are summed once per iteration (Section~\ref{sec:exp-request-scalability} reports iteration times). At the final temperature the earlier, unfused implementation, which makes about 20 separate sigmoid and scatter-reduce passes per batch, takes 4.76, 2.49, and 2.42~s per iteration on one, two, and four H100s, against 1.48, 0.75, and 0.39~s for the fused kernels, which are 3.2, 3.3, and 6.2 times faster. Both compute the same smoothed allocation: on the headline day their final bids differ by a median of $1.3\times10^{-6}$ and at most $9.9\times10^{-5}$ over active campaigns (Appendix~\ref{sec:kernel_correctness}).

\subsection{A log-bid update with momentum and annealing}\label{sec:bid_update}

\paragraph{Log-bid step.} The solver works in log-bid coordinates $x_j=\log b_j$. Each iteration makes one auction pass at the current bids $\mathbf x^{(t)}$, which gives the spends $S_j^{(t)}$, and forms the direction
\begin{equation}\label{eq:direction}
    g_j^{(t)}\;=\;\log\mathrm{clip}\!\left(\frac{B_j}{S_j^{(t)}},\ r_{\mathrm{lo}},\ r_{\mathrm{hi}}\right),\qquad [r_{\mathrm{lo}},r_{\mathrm{hi}}]=[0.1,10].
\end{equation}
Near equilibrium $\log(B_j/S_j)=-(S_j-B_j)/B_j+O\bigl((S_j/B_j-1)^2\bigr)$, so a step along $\mathbf g^{(t)}$ is, to first order, the projected gradient step of Algorithm~\ref{alg:solver} preconditioned by inverse budget, which makes the step dimensionless and puts small and large campaigns on the same scale. It is also a spending-ratio variant of the budget dynamics of \citet{borgs2007dynamics}, which raise or lower each multiplier by a fixed factor according to the sign of $S_j-B_j$. The ratio clamp bounds the step while spend is far from budget, in particular at the unpaced start.

\paragraph{Momentum.} The update adds heavy-ball momentum~\citep{polyak1964some} with the momentum coefficients of FISTA~\citep{beck2009fista}:
\begin{equation}\label{eq:update}
    \mathbf y^{(t)}=\mathbf x^{(t)}+\gamma_t\bigl(\mathbf x^{(t)}-\mathbf x^{(t-1)}\bigr),\qquad
    \mathbf x^{(t+1)}=\Pi_{\mathcal X}\bigl(\mathbf y^{(t)}+\alpha_{\text{step}}(t)\,\mathbf g^{(t)}\bigr),
\end{equation}
where $\gamma_t=(\nu_{t-1}-1)/\nu_t$ with $\nu_t=\bigl(1+\sqrt{1+4\nu_{t-1}^2}\bigr)/2$ and $\nu=1$ at every (re)start, and $\Pi_{\mathcal X}$ clips each coordinate to the bid box $[\log\varepsilon,\log\cp_j]$ of Section~\ref{sec:kkt}. The direction is evaluated at the current iterate $\mathbf x^{(t)}$, where the auction pass is made, not at the extrapolated point $\mathbf y^{(t)}$ as in FISTA, so the method is projected heavy-ball momentum with FISTA's momentum schedule and costs one auction pass per iteration. Momentum is off for the first three iterations, is reset when the anneal reaches its final temperature, and is reset whenever the momentum and the step point in opposite directions, $(\mathbf y^{(t)}-\mathbf x^{(t+1)})^\top(\mathbf x^{(t+1)}-\mathbf x^{(t)})>0$, the adaptive restart rule of \citet{odonoghue2015restart}.

\paragraph{Step-size coupling.} By Proposition~\ref{prop:convex_smooth}(iii), the gradient Lipschitz constant of $\Lsm$ on any compact convex set scales as $L_\mathbf T=O(1/T)$, which suggests a gradient step $\sim 1/L_\mathbf T\sim T$. Since $T_i=\alpha_{\text{temp}}(t)\cdot\mathrm{std}(\mathbf s_i)$, we set $\alpha_{\text{step}}(t)=c\cdot\alpha_{\text{temp}}(t)$ with $c=0.42$, the same for every campaign. The update \eqref{eq:update} uses momentum and a clipped log-ratio direction, so this is guidance rather than a theorem.

\paragraph{Annealing schedule.} The temperature coefficient decays geometrically,
\begin{equation}
    \alpha_{\text{temp}}(t)\;=\;\alpha_0\,(\alpha_{\text{end}}/\alpha_0)^{\min(t,t_{\text{horizon}})/t_{\text{horizon}}},
\end{equation}
from $\alpha_0=1.25$ to $\alpha_{\text{end}}=0.01$ over $t_{\text{horizon}}=80$ iterations, so equal numbers of iterations shrink $T$ by equal ratios. After the anneal, $T$ stays at its final value while the iterates continue to converge. Section~\ref{sec:exp-ablations} shows that momentum, annealing, and the choice of $c$ each matter.

\subsection{What the theory covers}\label{sec:scope}

Sections~\ref{sec:model} and~\ref{sec:fd_relaxation} concern the \emph{ideal} problem: fixed temperatures $T_i>0$, the exact soft allocation \eqref{eq:primal}, and the exact minimizer $\mathbf x^\star(\mathbf T)$ of $\Lsm$. The implementation departs from it in three ways.

\begin{enumerate}
    \item[(S1)] \textbf{Bid-dependent temperature.} The solver recomputes each temperature from the current scores, $T_i(\mathbf x,t)=\alpha_{\text{temp}}(t)\cdot\mathrm{std}(\mathbf s_i(\mathbf x))$, so differentiating $V_i(\mathbf s_i(\mathbf x);T_i(\mathbf x))$ picks up a term through $\partial T_i/\partial x_j$ and $\mathbf S-\mathbf B$ is no longer exactly the gradient. (Fixing $T_i$ from the initial scores, or from any bid-independent scale, would restore it.) The limit is nevertheless covered: after the anneal $\alpha_{\text{temp}}$ is constant, so if the iterates converge to $\bar{\mathbf x}$, the temperatures $\bar T_i:=T_i(\bar{\mathbf x},t)$ are fixed positive numbers (the kernel floors them at $10^{-6}$). At $\bar{\mathbf x}$ the update \eqref{eq:update} leaves every coordinate unchanged, so $g_j(\bar{\mathbf x})$ is zero inside the box, nonnegative at the cap, and nonpositive at the floor; since the ratio clamp preserves the sign of $\log(B_j/S_j)$, these are the KKT conditions of Section~\ref{sec:kkt} for $\Lsm$ at $\bar{\mathbf T}$. By convexity $\bar{\mathbf x}$ minimizes $\widetilde{\mathcal L}_{\bar{\mathbf T}}$ over $\mathcal X$, so $\bar{\mathbf x}=\mathbf x^\star(\bar{\mathbf T})$ once $\varepsilon\le\min_j\underline b_j$ (Lemma~\ref{lem:floor}), and Theorem~\ref{thm:convergence} applies with $\mathbf T=\bar{\mathbf T}$: the bid-dependent temperature changes the solver's path, not its limit.
    \item[(S2)] \textbf{Numerical stabilization.} The kernel clamps logits at $\pm50$, sets $p_{ij}=1$ when $K_i=M_i$ (exact by Lemma~\ref{lem:fd-dual}), and renormalizes $p\leftarrow\mathrm{clamp}(p\cdot K_i/\sum p,0,1)$. In the audit of Appendix~\ref{sec:kernel_correctness} these perturb capacity by at most $2.1\times10^{-5}$, and the renormalization factor lies in $[0.99998,1.00002]$ on the 9.1\% of requests where it is active.
    \item[(S3)] \textbf{Annealed continuation.} The temperature changes during the first iterations, one update per temperature, so the iterates do not minimize any single $\Lsm$ until the anneal ends. We do not prove that they converge; if they do, (S1) identifies the limit, and in practice the solver stops at a point that meets these conditions approximately (Section~\ref{sec:convergence_table}).
\end{enumerate}

Theorem~\ref{thm:convergence} therefore bounds the distance from the exact smoothed equilibrium to the hard one, not from the solver's output. Corollary~\ref{cor:certificate} closes this gap after the fact: it bounds the distance from the solver's actual output to the hard FPPE, whatever (S1)--(S3) do. The anchor we report is that output, a smoothed proxy for the hard equilibrium whose distance from it we certify.

\section{Experiments}\label{sec:empirical}

We ask five questions: how close the smoothed equilibrium and the solver's output are to the exact equilibrium, on synthetic markets small enough to solve exactly (Section~\ref{sec:exp-synthetic}); how fast and how far the solver converges on a production workload, and how far its output is from the hard equilibrium (Section~\ref{sec:convergence_table}); how it compares with other methods, including the same update on the hard top-$K$ allocation (Section~\ref{sec:exp-baselines}); which of its components matter (Section~\ref{sec:exp-ablations}); and how its cost grows with request volume and the number of GPUs (Section~\ref{sec:exp-request-scalability}).

\subsection{Setup and metrics}\label{sec:exp-setup}
Each production workload is built from one day of LinkedIn Feed traffic, in which per-request slot counts $K_i$ vary widely. Requests are resampled uniformly at an undisclosed rate, which may subsample requests or duplicate them, so the size of a workload does not reveal production volumes; we likewise report no absolute request or candidate counts and no spend. The headline run (Section~\ref{sec:convergence_table}) uses one day; the same section repeats the main configuration on 23 consecutive days built the same way, the comparison with other methods (Section~\ref{sec:exp-baselines}) and the ablations (Section~\ref{sec:exp-ablations}) use three of them, and Section~\ref{sec:exp-request-scalability} uses nested samples of the headline day.

\paragraph{Main configuration.} The update of Section~\ref{sec:bid_update} with its stated settings (anneal from $1.25$ to $0.01$ over 80 iterations, $c=0.42$, ratio clamp $[0.1,10]$; Table~\ref{tab:hyper}), for at most 800 iterations. Auction passes run in the fused kernel of Section~\ref{sec:solver}, with requests sharded across four NVIDIA H100 GPUs for the headline run and two for the multi-day and ablation runs. We report iteration time excluding input loading, which is dominated by I/O.

\paragraph{Metrics.} Let $\mathcal{A}:=\{j : B_j>0,\ S_j>0,\ b_j<\cp_j\}$ denote the set of \emph{active} campaigns at the current iteration: campaigns with a positive budget that won some mass and are not pinned at their bid cap. Define the relative spend deviation $\delta_j := |B_j/S_j - 1|$. We report:
\begin{itemize}
    \item $w_p$ := percentage of $j\in\mathcal{A}$ with $\delta_j< p/100$, the share of active campaigns within $p\%$ of budget; we report $w_1$, $w_{0.1}$, and the looser $w_{10}$.
    \item The median, 90th, 99th and 99.9th percentiles of $\{\delta_j : j\in\mathcal{A}\}$, writing $Q_{0.99}$ for the 99th; the largest deviation $\max_{j\in\mathcal{A}}\delta_j$; and $L_1$ and $L_2$, the mean and root mean square of the deviations.
    \item $N_{\mathrm{over}}$ := $|\{j:x_j=\log\cp_j \text{ and } S_j>B_j\}|$, the number of campaigns pinned at the bid cap while still overspending. Such a campaign violates the KKT conditions of \S\ref{sec:kkt} (its bid should move below the cap), so $N_{\mathrm{over}}>0$ means the run has not converged.
\end{itemize}
Zero-budget campaigns have no target and are excluded from all metrics. Capped campaigns are excluded from $w_p$ and the deviation statistics because, by the KKT conditions of \S\ref{sec:kkt}, they may legitimately underspend; overspending capped campaigns are counted by $N_{\mathrm{over}}$. Thus $w_p$ measures convergence \emph{among reachable, uncapped campaigns}, not over the full population.

\paragraph{Stopping rule.} A run meets the \emph{practical stopping rule} once the anneal has reached its final temperature and, for 10 consecutive iterations, $w_{1}\ge 99.99\%$, $L_1<10^{-4}$, $L_2<10^{-3}$, $Q_{0.99}<10^{-3}$, and $N_{\mathrm{over}}=0$.

\subsection{How close is the smoothed equilibrium to the exact one?}\label{sec:exp-synthetic}

Production markets are too large to solve the hard equilibrium exactly, so we test the smoothing on synthetic markets where we can.

\paragraph{Setup.} Each of 10 random instances has 400 requests and 30 campaigns. A request has 4--14 eligible campaigns and $K_i\in\{1,\dots,4\}$ slots with $K_i<M_i$; key-result rates are log-normal with a per-campaign factor, bid caps are log-normal, and each budget is a uniform random fraction, between 0.15 and 1.3, of the campaign's spend when every campaign bids its cap, so some campaigns are budget-constrained and others end at their cap. We solve the hard problem of Theorem~\ref{thm:hard-dual} as a conic program with interior-point and splitting solvers~\citep{goulart2024clarabel,odonoghue2016scs}, to a KKT residual of at most $5\times10^{-7}$ on every instance. For fixed temperatures $T_i=\alpha\,\tau_i$, with $\tau_i$ the standard deviation of request $i$'s scores at the caps, we minimize the smoothed dual $\Lsm$ with L-BFGS-B~\citep{byrd1995lbfgsb} and its exact gradient $\mathbf S-\mathbf B$. Finally we run the practical solver of Section~\ref{sec:bid_update} from the caps for 800 iterations, with the main configuration's constants (Table~\ref{tab:hyper}) and temperatures recomputed from the current scores as in (S1). The experiment uses no production data.

\begin{figure}[!htbp]
\centering
\includegraphics[width=\textwidth]{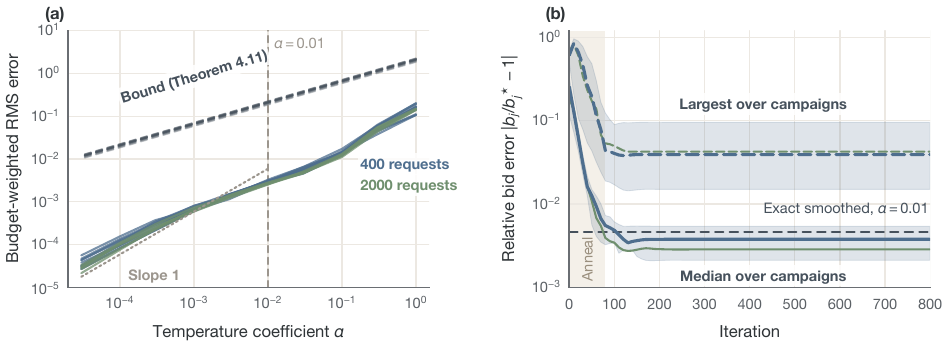}
\caption{Synthetic markets with an exactly computed FPPE. (a) Budget-weighted RMS error of the pacing multipliers of the exact smoothed equilibrium, one line per instance. (b) Relative bid error of the solver against the exact FPPE (median over instances; band: range over the 400-request instances).}
\label{fig:synthetic}
\end{figure}

\paragraph{Results.} Figure~\ref{fig:synthetic} summarizes the results.
\begin{itemize}
    \item \emph{The bound holds with a wide margin.} The budget-weighted error $\bigl(\sum_jB_j(\beta_j^\star(\mathbf T)-\beta_j^\star)^2\bigr)^{1/2}$ is at most 9.3\% of the bound $2\sqrt{\Gamma_\mathbf T}$ of Theorem~\ref{thm:convergence} at every temperature, and at most 2\% for $\alpha\le0.03$.
    \item \emph{The error decays linearly in $T$.} Over $\alpha\in[3\times10^{-5},3\times10^{-4}]$, the median relative bid error $|b_j^\star(\mathbf T)/b_j^\star-1|$ has fitted log--log slopes between 0.96 and 1.27 (median 1.03), so on these instances the worst-case $O(\sqrt T)$ rate is conservative.
    \item \emph{The final temperature sets the anchor's accuracy.} At $\alpha=0.01$, the value our schedule ends at, the smoothed equilibrium bids differ from the exact ones by a median of 0.47\% (median over instances; worst campaign 12\%); at $\alpha=0.001$ the median is 0.07\% (worst 2.2\%).
    \item \emph{The solver reaches the smoothed equilibrium.} The practical solver's final bids differ from the exact FPPE bids by a median of 0.38\% (worst campaign 9.5\%), no worse than the exact smoothed equilibrium at $\alpha=0.01$, while its median spend deviation is below $10^{-14}$. On these instances the approximations (S1)--(S3) cost nothing measurable, and the remaining error is the smoothing bias of the final temperature.
\end{itemize}
Five larger instances (2000 requests, 100 campaigns), whose hard solves are less precise (KKT residuals up to $3\times10^{-5}$), give the same picture: the error stays below 7.2\% of the bound, the smoothed bids at $\alpha=0.01$ differ from the exact ones by a median of 0.35\%, and the practical solver's by 0.29\%.

Budget matching and anchor accuracy are thus different quantities: spend can match budgets to $10^{-6}$ while bids carry the $O(T)$ smoothing bias, which lowering the final temperature reduces.

\subsection{How far does the solver converge at production scale?}\label{sec:convergence_table}
Table~\ref{tab:headline} summarizes the main configuration on the headline day, run for 800 iterations with the stopping rule switched off so that its behavior after convergence is visible. The rule would fire after 174 iterations, about 61~s on four H100s excluding input loading; by then every active campaign is within 1\% of its budget, 99.99\% are within 0.1\%, and no capped campaign overspends. Along the way, $w_1$ exceeds 99\% at iteration 96 and the median deviation falls below $10^{-3}$ at iteration 109; at iteration 153 the median falls below $10^{-4}$ and every active campaign comes within 1\% of its budget, where it stays. Continuing to 800 iterations lowers the median and the 99th percentile, while the largest deviation rises. At the end, 17.1\% of budgeted campaigns sit at their bid cap, where by \S\ref{sec:kkt} they may underspend at an optimum. At the final bids, allocating every request by hard top-$K$ instead raises total spend by only 0.011\%: in aggregate the smoothed allocation is nearly sharp at the final temperature.

\begin{table}[htbp]
\centering
\small
\caption{Main configuration on the headline day, four H100 GPUs. Deviations are $\delta_j=|B_j/S_j-1|$ over active campaigns; iteration time excludes input loading.}
\label{tab:headline}
\begin{tabular}{lrr}
\toprule
Metric & At the stopping rule & After 800 iterations \\
\midrule
Iterations & 174 & 800 \\
Iteration time, mean (median; p99) & \multicolumn{2}{c}{0.38~s (0.39~s; 0.44~s)} \\
Total iteration time & 61~s & 306~s \\
$w_{10}$ / $w_{1}$ / $w_{0.1}$ & 100\% / 100\% / 99.994\% & 100\% / 100\% / 99.995\% \\
Median $\delta_j$ & $4.0\times10^{-5}$ & $8.8\times10^{-6}$ \\
90th percentile of $\delta_j$ & $1.2\times10^{-4}$ & $3.1\times10^{-5}$ \\
99th percentile of $\delta_j$ & $2.0\times10^{-4}$ & $1.1\times10^{-4}$ \\
99.9th percentile of $\delta_j$ & $3.4\times10^{-4}$ & $3.1\times10^{-4}$ \\
Largest $\delta_j$ & $2.8\times10^{-3}$ & $5.8\times10^{-3}$ \\
$L_1$ / $L_2$ of $\delta_j$ & $5.5\times10^{-5}$ / $7.5\times10^{-5}$ & $1.5\times10^{-5}$ / $3.2\times10^{-5}$ \\
$N_{\mathrm{over}}$ & 0 & 0 \\
\bottomrule
\end{tabular}
\end{table}

\begin{figure}[!htbp]
\centering
\includegraphics[width=\textwidth]{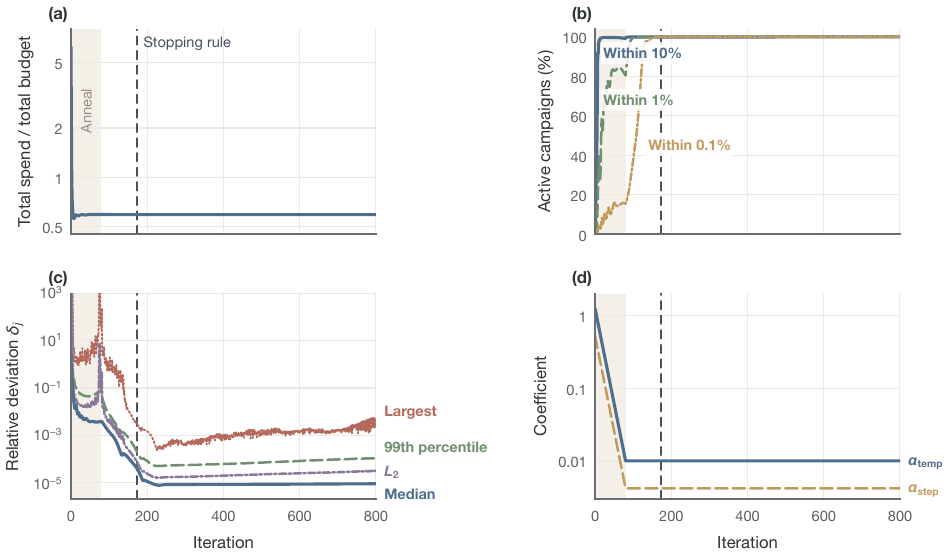}
\caption{Convergence on the headline day (main configuration). (a) Total spend over total budget. (b) Active campaigns within 10\%, 1\%, and 0.1\% of budget. (c) Relative deviation over active campaigns; axis cut at $10^3$. (d) Temperature and step coefficients.}
\label{fig:convergence}
\end{figure}

Figure~\ref{fig:convergence} shows three phases. During the anneal (iterations 0--80), the unpaced initial bids, which overspend about sixfold, are brought to budget within about ten iterations, and $w_{10}$ exceeds 99\% by iteration 14; the tighter bands fill only once the temperature is low. When the anneal ends, a few campaigns cause a brief spike in the RMS deviation, after which every metric falls steadily until the stopping rule fires. After that the run is flat: the median settles near $10^{-5}$, the tail creeps up slightly (from iteration 300 to 800, the 99th percentile rises from $5\times10^{-5}$ to $1.1\times10^{-4}$ and the largest deviation from $4.6\times10^{-4}$ to $5.8\times10^{-3}$), and every active campaign stays within 1\% of budget.

\paragraph{Twenty-three consecutive days.} We ran the main configuration on 23 consecutive days, on two H100s, for 800 iterations with the stopping rule switched off. On 14 days the rule would fire after 157--200 iterations (median 162). On the other 9, every criterion except $N_{\mathrm{over}}=0$ is met after 157--165 iterations and stays met to the end of the run: one to four campaigns sit at their cap overspending by at most $7.0\times10^{-6}$ of their budget, which the zero-tolerance criterion still counts. Ignoring that criterion, the rule fires after 157--185 iterations on every day. At that point $w_1>99.99\%$ and $w_{0.1}>99.9\%$ on every day, and on 18 of the 23 days every active campaign is within 1\% of its budget; by iteration 200 this holds on every day, with the largest deviation at most $3.2\times10^{-3}$. After 800 iterations, every active campaign is within 1\% of its budget and $w_{0.1}\ge99.99\%$ on every day, the median deviation is between $7.1\times10^{-6}$ and $1.1\times10^{-5}$, and the 99th percentile is at most $1.9\times10^{-4}$; no run diverges or oscillates. An iteration takes 0.50--1.12~s on two H100s, depending on the day's volume.

\paragraph{Distance to the hard equilibrium.} Corollary~\ref{cor:certificate} turns the final bids of each run into a bound on their distance from the hard FPPE. We evaluate it with one hard top-$K$ pass and smoothed passes at several temperatures at the final bids, taking for $\mathbf p$ the best convex combination of two of these allocations. On the headline day the gap $\Delta$ is $6.2\times10^{-5}$ of total budget, so the budget-weighted root-mean-square error of the pacing multipliers, $\bigl(\sum_jB_j(\beta_j-\beta^\star_j)^2/\sum_jB_j\bigr)^{1/2}$, is at most $0.016$, which is 2.2\% of their budget-weighted root-mean-square value. Across the 23 days the gap is between $4.7\times10^{-5}$ and $8.2\times10^{-5}$ of total budget and the bound between $0.014$ and $0.018$ (1.9--2.8\% relative). The gap almost equals the revenue difference between the hard and the smoothed allocations at the final bids ($6.4\times10^{-5}$ of total budget on the headline day), so the bound is set by the smoothing at the final temperature, not by the solver's residual; it is worst-case and aggregate, and individual anchors may lie further from or closer to the hard equilibrium. The kernels enforce capacity only up to single-precision rounding; deflating the smoothed allocation by the largest capacity residual of the audit in Appendix~\ref{sec:kernel_correctness}, $2.1\times10^{-5}$, raises every bound by at most $0.003$.

\paragraph{Stability of the anchors.} On three of these days we compared bids across runs. Starting from bids multiplied by independent random factors $e^{U[-1,1]}$, up to 2.7$\times$ in either direction, the final bids differ from the unperturbed run's by a median of at most $2.5\times10^{-5}$ (99th percentile at most $1.7\times10^{-4}$). Continuing to 1500 iterations changes them by a median of at most $2.3\times10^{-5}$. Repeated runs are bit-identical, and the headline day gives identical per-iteration metrics on one, two, and four GPUs (Section~\ref{sec:exp-request-scalability}). The anchors are therefore set by the smoothed problem, not by the solver's path, and their remaining error is the smoothing bias of the final temperature: lowering it from 0.01 to 0.001 moves the bids by a median of 0.1\% and a 99th percentile of 1.3\%, in line with the synthetic estimates of Section~\ref{sec:exp-synthetic}.

\subsection{How does the solver compare with other methods?}
\label{sec:exp-baselines}\label{sec:exp-smoothing-ablation}

\paragraph{Setup.}
On three of the 23 days, which we call Days A, B, and C (Day B is the headline day), we run each
method for 800 iterations from the same initial bids with the stopping rule switched off.
\begin{itemize}
    \item \emph{On the smoothed allocation:} the main configuration; Algorithm~\ref{alg:solver},
    projected gradient with the exact gradient $\mathbf S-\mathbf B$ and our temperature schedule, whose
    scalar step $0.42\,\alpha_{\text{temp}}(t)$ we divide by the median or by the largest budget, because
    $\mathbf S-\mathbf B$ is measured in units of spend; and Adam~\citep{kingma2014adam} with the same exact
    gradient in log-bid coordinates, at constant learning rates $0.01$, $0.03$, and $0.1$ and at our step
    schedule.
    \item \emph{On the hard top-$K$ allocation:} our update, unchanged, to isolate the effect of the
    allocation; and two methods of \citet{afiouni2026pacing}. \emph{PACE} is dual averaging on the hard
    dual, which they derive from \citet{gao2021online}: it sets $\beta_j\leftarrow\min\{B_j/\bar u_j,1\}$,
    where $\bar u_j$ is the mean over all past iterations of the value $\sum_iv_{ij}p_{ij}$ that campaign
    $j$ receives. Their version of the \emph{budget dynamics} of \citet{borgs2007dynamics} multiplies
    $\beta_j$ by $e^{-\eta}$ when campaign $j$ overspends and by $e^{\eta}$ otherwise, in an auction in
    which every bid is multiplied by an independent $U[0.8,1]$ factor at every iteration, the perturbation
    size of their experiments; we use $\eta\in\{0.01,0.03\}$.
\end{itemize}
The methods use different allocations, so budget residuals alone do not compare them.
Corollary~\ref{cor:certificate} does: for each run we evaluate the hard dual $\Lhard+C_0$ at its final
bids and subtract the largest primal value found by any run that day, which by the corollary is a valid
gap $\Delta$ for every run. We report $\Delta$ relative to total budget; $2\sqrt{\Delta/\sum_jB_j}$ bounds
the budget-weighted root-mean-square error of the pacing multipliers. The fused kernel computes only the
smoothed allocation, so the hard-allocation methods run on the earlier, unfused kernel, and we do not
compare iteration times.

\begin{table}[htbp]
\centering
\small
\setlength{\tabcolsep}{4pt}
\caption{After 800 iterations on three days. $w_1$: share of active campaigns within 1\% of budget under each method's own allocation (budget dynamics: one random draw). $\Delta/\sum_jB_j$: certificate of Corollary~\ref{cor:certificate} at the final bids.}
\label{tab:baselines}
\begin{tabular}{llcc}
\toprule
Method & Allocation & $w_1$ (\%), Days A / B / C & $\Delta/\sum_jB_j$, Days A / B / C \\
\midrule
Main configuration & smoothed & 100 / 100 / 100 & $7.0$ / $6.2$ / $7.6\times10^{-5}$ \\
Alg.~\ref{alg:solver}, median-budget step & smoothed & 27.9 / 25.4 / 25.9 & $0.93$ / $12.7$ / $3.6$ \\
Alg.~\ref{alg:solver}, largest-budget step & smoothed & 0.4 / 0.4 / 0.4 & $0.47$ / $0.40$ / $0.47$ \\
Adam, rate $0.1$ & smoothed & 81.0 / 80.5 / 80.6 & $3.5$ / $3.9$ / $4.2\times10^{-4}$ \\
Adam, rate $0.03$ & smoothed & 97.0 / 96.3 / 96.9 & $5.5$ / $6.8$ / $8.0\times10^{-3}$ \\
Adam, rate $0.01$ & smoothed & 94.0 / 91.6 / 94.3 & $1.1$ / $1.3$ / $1.4\times10^{-1}$ \\
Adam, our step schedule & smoothed & 99.0 / 98.8 / 99.0 & $3.6$ / $4.4$ / $5.3\times10^{-3}$ \\
Our update & hard & 51.4 / 44.0 / 50.1 & $1.5$ / $1.4$ / $1.7\times10^{-3}$ \\
PACE & hard & 73.3 / 68.3 / 72.3 & $6.6$ / $6.7$ / $7.8\times10^{-5}$ \\
Budget dynamics, $\eta=0.01$ & perturbed hard & 2.2 / 2.3 / 2.2 & $1.6$ / $1.3$ / $1.7\times10^{-2}$ \\
Budget dynamics, $\eta=0.03$ & perturbed hard & 2.2 / 2.2 / 2.2 & $2.2$ / $1.7$ / $2.3\times10^{-2}$ \\
\bottomrule
\end{tabular}
\end{table}

\begin{figure}[!t]
    \centering
    \includegraphics[width=\linewidth]{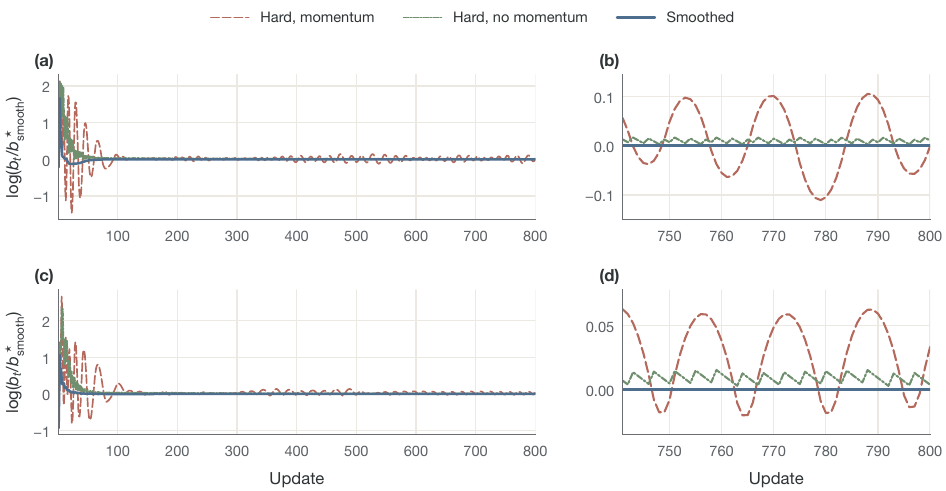}
    \caption{Bids on Day B of the two traced campaigns (rows) that vary most under hard top-$K$ without momentum, normalized by their final smoothed bids: (a, c) all 800 updates; (b, d) the last 60.}
    \label{fig:exp-ablation-convergence}
\end{figure}

\paragraph{Smoothing.}
With the update fixed, the allocation decides the outcome. On Day B, smoothing brings every active
campaign within 1\% of its budget with a median deviation of $8.8\times10^{-6}$, the same as the fused
run of Section~\ref{sec:convergence_table}. Under hard top-$K$, only 44.0\% of active campaigns end within
1\% of budget and 84.3\% within 10\%, the median deviation is $1.3\times10^{-2}$, 84 capped campaigns still
overspend, and bids are still moving by large amounts at the last update; the certificate is about twenty
times the main configuration's on every day (Table~\ref{tab:baselines}). The failure is not caused by
momentum: without it, hard top-$K$ leaves 67.5\% of active campaigns within 1\% of budget on Day B and 93
capped campaigns overspending, whereas the smoothed update without momentum keeps 99.99\% within 1\% but
stalls (Section~\ref{sec:exp-ablations}). Smoothing therefore helps with either update, while momentum
helps only when the allocation is smooth.

The hard-allocation runs update all campaigns jointly but record the bids of only a random 0.1\% of
them. Figure~\ref{fig:exp-ablation-convergence} shows the two traced campaigns whose bids vary most over
the last 300 updates of hard top-$K$ without momentum, under hard top-$K$ with and without momentum and
under smoothing with momentum. Without momentum, hard top-$K$ produces
a sawtooth: the bid rises by one large step when the campaign loses its marginal auctions and
underspends, then falls in several small steps while it overspends, and the cycle repeats every few
updates. Momentum turns the sawtooth into a slower oscillation of larger amplitude. The smoothed bids
settle. Among the 259 traced campaigns that end active, the standard deviation of the log bid over the
last 300 updates has a median of $2.6\times10^{-3}$ and a maximum of $0.22$ under hard top-$K$ with
momentum and a median of $9.3\times10^{-4}$ and a maximum of $4.5\times10^{-3}$ without, against at most
$9.4\times10^{-6}$ with smoothing.

\paragraph{Other methods.}
The main configuration is the only method in Table~\ref{tab:baselines} that brings every active
campaign within 1\% of its budget, and its certificate is the smallest on Days~B
and~C. PACE comes close on the hard dual: its certificate is within 10\% of ours on all three days and
6\% smaller on Day~A, and its final bids differ from ours by a median of 0.2\% (99th percentile 1.5\%),
the order of the smoothing bias at our final temperature. Its iterates, however, are still moving after
800 iterations, with a median budget deviation of $3\times10^{-3}$ that is still falling, and under the
hard allocation only 68--73\% of active campaigns end within 1\% of budget: dual averaging finds nearly
the same bids but not an allocation that clears the budgets at them, which on a near-tie requires
splitting slots fractionally. Algorithm~\ref{alg:solver}, the reference method of
Section~\ref{sec:reference}, fails in both normalizations, because a single step cannot suit budgets that
span orders of magnitude: scaled by the median budget, the largest campaigns take steps far too large,
oscillate, and leave total spend 5--16 times total budget; scaled by the largest budget, most campaigns
barely move and fewer than 0.5\% come within 1\% of budget. Dividing the gradient by each campaign's own
budget, as our update does to first order, removes this mismatch. Adam, given the same exact gradient,
meets the stopping rule at none of the four learning rates: the best certificate, at rate $0.1$, is five to
six times ours, and that run leaves a fifth of active campaigns more than 1\% from budget. The budget
dynamics settle where the perturbed auction clears budgets in expectation, so their bids differ from ours
by a median of 17--18\%, and their certificate is two to three hundred times ours.

\subsection{Which solver components matter?}\label{sec:exp-ablations}

\paragraph{Setup.}
Starting from the main configuration, we change one component at a time and run each variant on Days A--C of Section~\ref{sec:exp-baselines}, on two H100s with the fused kernel. Each run takes 800 iterations with the stopping rule switched off, and Table~\ref{tab:ablations} reports the iteration at which the practical stopping rule of Section~\ref{sec:exp-setup} would fire. Each cell is a single run; the solver is deterministic, so a repeated run gives the same result.

\begin{table}[htbp]
\centering
\small
\setlength{\tabcolsep}{4pt}
\caption{Iterations until the stopping rule fires (at most 800), changing one component of the main configuration at a time. $^\dagger$All criteria met except $N_{\mathrm{over}}=0$ (one capped campaign overspends marginally).}
\label{tab:ablations}
\begin{tabular}{lcccl}
\toprule
Variant & Day A & Day B & Day C & Notes \\
\midrule
Main configuration ($c=0.42$) & \textbf{170} & \textbf{174} & \textbf{171} & \\
$c=0.30$ & $>800^\dagger$ & $>800^\dagger$ & 206 & \\
$c=0.60$ & $>800^\dagger$ & 197 & 204 & a few campaigns keep oscillating \\
$c=0.90$ & --- & --- & --- & campaigns oscillate ($w_1$ 98.8--99.95\% at 800) \\
$c=1.35$ & --- & --- & --- & diverges ($w_1$ 36--45\% at 800) \\
Anneal length 40 & $>800^\dagger$ & 163 & 171 & \\
Anneal length 160 & $>800^\dagger$ & 222 & 228 & \\
Final temperature $\alpha_{\text{end}}=0.005$ & 211 & 213 & 196 & \\
No momentum & --- & --- & --- & stalls: median $\delta_j\approx5\times10^{-4}$ \\
No anneal ($\alpha_{\text{temp}}\equiv0.01$) & --- & 430 & --- & unstable ($w_1$ 68\% at 800) \\
Constant $\alpha_{\text{temp}}\equiv0.05$ & 535 & 115 & 429 & final temperature 5$\times$ higher \\
\bottomrule
\end{tabular}
\end{table}

\paragraph{Results.}
Momentum and annealing are both essential. Without momentum the solver stalls at a median deviation of about $5\times10^{-4}$ and never meets the rule; on Day B that is fifty times the main configuration's median after the same 800 iterations, even though its final bids differ from the main configuration's by a median of only 0.04--0.05\% and its certificate is within 3\% of the main configuration's on every day. Starting directly at the final temperature is unstable on two of the three days; a constant, higher temperature converges, but to a softer problem whose smoothing bias is about five times larger (Section~\ref{sec:exp-synthetic}). The step coefficient has a narrow good range: $c=0.42$ stops first on every day and smaller steps are slower. Larger steps destabilize a few campaigns: with $c=0.60$ the rule fires on two days, but one to nine campaigns stay outside the 1\% band, oscillating, for the rest of the run; with $c=0.90$ up to about 1\% of active campaigns never settle; and three times the step diverges. A shorter anneal gains nothing and a longer one costs about 50 iterations. Halving the final temperature, which roughly halves the smoothing bias, costs 20--40 iterations. Several variants never stop only because one capped campaign keeps overspending marginally; Section~\ref{sec:convergence_table} shows that this also happens with the main configuration on some days and is harmless.

\subsection{How does the cost grow with request volume?}
\label{sec:exp-request-scalability}

\paragraph{Setup.}
From the headline day we draw nested uniform samples of requests: each request
is kept independently with probability $1/64$, $1/16$, $1/4$, or $1/2$, with one
draw per request shared across rates, so that each sample contains the smaller
ones. A kept request retains its complete candidate list and slot count; bids
and caps are unchanged, and every budget is multiplied by the sampling rate, so
each campaign faces about the same budget pressure as on the full day. On
each sample and on the full day we run the main configuration for 800
iterations, with the stopping rule switched off, on one, two, and four H100
GPUs, twice for each combination, and report the mean time of a complete
iteration, which excludes input loading. The two runs of a combination differ
by at most 3.2\%.

\begin{figure}[H]
    \centering
    \includegraphics[width=\linewidth]{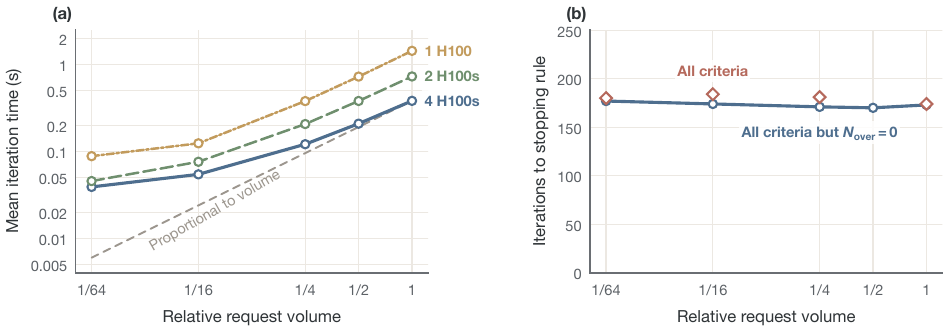}
    \caption{Cost against request volume on nested samples of the headline day. (a) Mean iteration time, excluding loading (two runs each). (b) Iteration at which the stopping rule fires; the same on every GPU count.}
    \label{fig:exp-request-scalability}
\end{figure}

\paragraph{Results.}
The iteration time grows in proportion to the request volume once the GPUs are
busy, and falls in proportion to their number
(Figure~\ref{fig:exp-request-scalability}a). On the full day an iteration
takes 1.44~s on one H100, 0.73~s on two, and 0.38~s on four. From a quarter of
the day up, doubling the volume multiplies the iteration time by 1.73--1.98, and
halving the number of GPUs multiplies it by 1.71--1.97; both ratios approach two
as the volume grows. On the smallest sample a fixed cost per iteration
dominates: at $1/64$ of the day an iteration takes 0.09~s on one H100 and
0.04--0.05~s on two or four. The number of iterations does not grow with the volume
(Figure~\ref{fig:exp-request-scalability}b): every criterion of the stopping rule except $N_{\mathrm{over}}=0$ is
met after 170--177 iterations at every rate, and the full rule fires after
174--184 iterations at every rate except $1/2$, where one capped campaign keeps
overspending marginally. After 800 iterations every active campaign is within
1\% of its budget at every rate. The time to solve is therefore roughly
proportional to the request volume divided by the number of GPUs. At every
rate, the runs on one, two, and four GPUs give identical per-iteration metrics.

\section{Discussion}\label{sec:discussion}

\paragraph{Limitations.} The solver departs from the exact smoothed problem in the ways (S1)--(S3) of Section~\ref{sec:scope}, so the anchor is a smoothed proxy for the hard equilibrium. Corollary~\ref{cor:certificate} bounds its distance from it after the fact, but only in aggregate: at production scale the budget-weighted error of the pacing multipliers is at most about 2\% of their size (Section~\ref{sec:convergence_table}), and the bound says nothing about any single campaign. Campaign by campaign, on synthetic markets (Section~\ref{sec:exp-synthetic}), the solver reaches the smoothed equilibrium, whose bids lie within about 0.5\% of the exact ones at our final temperature, in line with how the production anchors move as the final temperature falls (Section~\ref{sec:convergence_table}). The production evidence covers 23 consecutive days of one feed-ad workload; the solver is deterministic, so repeated runs agree exactly, but we have not studied other surfaces. The anchor also depends on the traffic it is computed from: on random subsamples of the workload, the anchors of campaigns with ample traffic are stable to about a percent, while those of campaigns with little traffic are noisy. The model prices a request's slots first-price, treats them as interchangeable, and takes the platform's value predictions as given; these are modeling choices for tractability, not claims about any deployed auction (Section~\ref{sec:why_fpa}). The anchor is the equilibrium of the model market, a counterfactual reference rather than a forecast of production outcomes, and it inherits any miscalibration of the predicted click-through rates. We do not evaluate its downstream uses (Section~\ref{sec:intro}). Finally, the anchor is informative only where the budget binds: in the main run 17.1\% of budgeted campaigns end at their bid cap, where the anchor is the cap itself.

\paragraph{Faster solvers.} For a fixed temperature the exact gradient \eqref{eq:SminusB} makes classical convex optimization directly applicable. Limited-memory BFGS with box constraints (L-BFGS-B)~\citep{byrd1995lbfgsb} is a natural candidate: each forward pass returns both $\Lsm$ and $\mathbf S-\mathbf B$, so it obtains genuine secant pairs and learns curvature without a finite-difference estimator. With adaptive continuation in $\mathbf T$ it could replace the ratio clamp and the hand-set step schedule; the narrow range of good step sizes in Section~\ref{sec:exp-ablations} suggests room to gain, although Adam with the same exact gradient does worse than our update (Section~\ref{sec:exp-baselines}). PACE dual averaging on the hard allocation comes close to our certificate without clearing budgets (Section~\ref{sec:exp-baselines}); averaging over iterations on the smoothed allocation might combine the two. Lowering the final temperature would reduce the anchor's smoothing bias, linear in $T$ on synthetic markets (Section~\ref{sec:exp-synthetic}), and sharpen the allocations (at our final temperature, 18.2\% of audited requests still place more than $0.01$ of their mass outside the top $K_i$; Appendix~\ref{sec:kernel_correctness}); it would likely call for Newton steps, which are affordable because the Hessian \eqref{eq:hessian-sandwich} is a diagonal plus one rank-one term per request, so a Hessian--vector product costs one extra pass. Accuracy is ultimately limited by the kernel's \texttt{float32} arithmetic, roughly $10^{-7}$ in relative spend.

\paragraph{Open theoretical questions.} (i)~\emph{Position-dependent values:} allowing slot-specific click-through rates replaces the capped simplex by an assignment polytope, where \citet{afiouni2026pacing} show that the market-equilibrium interpretation can fail. (ii)~\emph{Return-on-spend constraints:} autobidders often carry ROI targets alongside budgets~\citep{aggarwal2019autobidding,deng2024efficiency}. (iii)~\emph{Online multi-slot pacing:} the smoothed gradient could drive dual-descent pacing online~\citep{balseiro2023best}. (iv)~\emph{Sharper rates:} the $O(\sqrt{T})$ bound of Theorem~\ref{thm:convergence} is worst-case, while the synthetic error decays linearly in $T$ (Section~\ref{sec:exp-synthetic}); a bound that exploits non-degenerate boundary gaps, as in Lemma~\ref{lem:mismass}, may explain the linear rate.

\section{Conclusion}\label{sec:conclusion}

The multi-slot first-price pacing equilibrium is unique and solves a convex program over bids. We showed that it can also be computed at production scale by smoothing each request's allocation with the binary entropy: the smoothing keeps the dual gradient exactly equal to spend minus budget, reduces each iteration to independent per-request threshold solves that parallelize on a GPU, and converges to the hard equilibrium at an explicit $O(\sqrt{T})$ rate. A solver built on this structure meets our stopping rule on a production feed-ad workload after 174 iterations, about a minute on four H100 GPUs, with a median relative spend deviation of $4.0\times10^{-5}$, and reaches the same accuracy within 200 iterations on each of 23 consecutive days; with the hard allocation and the same update, 56.0\% of active campaigns remain more than 1\% from budget and bids keep cycling. A duality-gap certificate bounds the budget-weighted error of the computed pacing multipliers by under 3\% of their size, and of the other methods we compare only PACE dual averaging comes close, without clearing budgets. The resulting equilibrium price anchor gives each campaign a hindsight reference for its cost per key result, against which online pacing, unit economics, and proposed market changes can be measured.

\bibliographystyle{plainnat}
\bibliography{references}

\appendix

\section{Deferred Proofs for Section~\ref{sec:model}}\label{app:proofs-hard}

\subsection{Proof of Proposition~\ref{prop:slater}}\label{app:slater}

\paragraph{Feasibility.} Set $p_{ij}=0$, $\omega_j=B_j$, and $u_j=B_j/2$. All constraints hold and $u_j>0$.

\paragraph{Strong duality and dual attainment.} Program \eqref{eq:pk} maximizes a concave objective whose domain $\{u_j>0\}$ is open, subject to affine constraints. The point above is feasible and lies in the interior of the objective's domain, so the refined Slater condition holds~\citep[Section~5.2.3]{boyd2004convex}: whenever the optimal value is finite, there is zero duality gap and the dual optimum is attained. Finiteness follows from the attainment argument below. The dual function is computed in Section~\ref{ssec:lagrangian}: the supremum of \eqref{eq:lagrangian} over $u_j>0$ equals $B_j\log(B_j/\beta_j)-B_j$ if $\beta_j>0$ and $+\infty$ otherwise; over $\omega_j\ge 0$ it is $0$ if $\beta_j\le 1$ and $+\infty$ otherwise; and over $p_{ij}\in[0,1]$ it is $\hinge{\beta_jv_{ij}-c_i}$.

\paragraph{Primal attainment.} Let $\Psi(\mathbf p,\mathbf u,\bm\omega):=\sum_jB_j\log u_j-\sum_j\omega_j$, let $\underline\Psi$ be its value at the feasible point above, and let $V_j^{\max}:=\sum_iv_{ij}$. Any maximizing sequence may be restricted to the feasible superlevel set $\{\Psi\ge\underline\Psi\}$. For feasible points, $u_j\le V_j^{\max}+\omega_j$, so $\Psi\le\sum_j\varphi_j(\omega_j)$ with $\varphi_j(t):=B_j\log(V_j^{\max}+t)-t$. Each $\varphi_j$ tends to $-\infty$ as $t\to\infty$, so $\bar\varphi_j:=\sup_{t\ge 0}\varphi_j(t)<\infty$. On the superlevel set, $\varphi_j(\omega_j)\ge\underline\Psi-\sum_{k\ne j}\bar\varphi_k$ for each $j$, which forces $\omega_j\le D_j$ for a finite constant $D_j$; hence $u_j\le V_j^{\max}+D_j$. Also $\underline\Psi\le\Psi\le B_j\log u_j+\sum_{k\ne j}B_k\log(V_k^{\max}+D_k)$, so $u_j\ge m_j$ for a constant $m_j>0$. With $0\le p_{ij}\le 1$, the feasible part of the superlevel set is closed and bounded, hence compact, and $\Psi$ is continuous on it. The maximum is attained, and in particular the optimal value is finite. \hfill$\square$

\subsection{Proof of Lemma~\ref{lem:topk}}\label{app:topk}

Let $\phi(c):=Kc+\sum_j\hinge{s_j-c}$.

\paragraph{Lower bound.} For any $c\in\R$, using $\hinge{z}\ge z$ on the top $K$ terms and $\hinge{z}\ge 0$ on the rest,
\begin{equation*}
    \phi(c)\ge Kc+\sum_{m=1}^K(s_{(m)}-c)=\sum_{m=1}^Ks_{(m)}.
\end{equation*}

\paragraph{Attainment.} For $c\in[s_{(K+1)},s_{(K)}]$, $\hinge{s_{(m)}-c}=s_{(m)}-c$ for $m\le K$ and $0$ for $m>K$, so $\phi(c)=\sum_{m=1}^Ks_{(m)}$. Since $s_{(K+1)}\ge 0$, these thresholds are feasible for $c\ge 0$.

\paragraph{No other minimizers.} If $0\le c<s_{(K+1)}$, at least $K+1$ scores exceed $c$ and $\phi(c)\ge\sum_{m=1}^{K+1}(s_{(m)}-c)+Kc=\sum_{m=1}^Ks_{(m)}+(s_{(K+1)}-c)$, which is strictly larger. If $c>s_{(K)}$, fewer than $K$ scores exceed $c$; writing $a<K$ for their number, $\phi(c)=\sum_{m\le a}s_{(m)}+(K-a)c>\sum_{m=1}^Ks_{(m)}$ because $s_{(m)}\le c$ for $a<m\le K$, with strict inequality at $m=K$.

\paragraph{Minimization over $\R$.} If $K<n$ and $c<0$, then $\phi(c)=\sum_js_j-(n-K)c>\sum_js_j\ge\TopK_K(\mathbf s)$, so negative thresholds are never optimal.

\paragraph{Properties.} $\TopK_K(\mathbf s)=\max_{|A|=K}\sum_{j\in A}s_j$ is a maximum of linear functions with non-negative coefficients, hence convex, positively homogeneous, and non-decreasing in each coordinate. \hfill$\square$

\subsection{Proof of Theorem~\ref{thm:hard-dual}}\label{app:hard-dual}
By Lemma~\ref{lem:topk}, minimizing \eqref{eq:dk} over $\mathbf c$ for fixed $\bm\beta$ gives $\sum_i\TopK_{K_i}(\bm\beta\odot\mathbf v_i)-\sum_jB_j\log\beta_j+C_B$; changing variables to $x_j=\log(\beta_j\cp_j)$ gives $\Lhard(\mathbf x)+C_0$. Attainment and the equality follow from Proposition~\ref{prop:slater}. Convexity: $\TopK_{K_i}$ is convex and non-decreasing in each coordinate and each $x_j\mapsto r_{ij}e^{x_j}$ is convex, so each term is convex; the remaining term is linear. \hfill$\square$

\subsection{Proof of Theorem~\ref{thm:characterization}}\label{app:characterization}

\paragraph{Optimality implies equilibrium.} Let $(\mathbf p,\mathbf u,\bm\omega)$ be optimal for \eqref{eq:pk} and $(\bm\beta,\mathbf c)$ optimal for \eqref{eq:dk}; both exist by Proposition~\ref{prop:slater}. By zero duality gap, $(\mathbf p,\mathbf u,\bm\omega)$ maximizes the Lagrangian \eqref{eq:lagrangian} at $(\bm\beta,\mathbf c)$, and complementary slackness holds. Set $b_j=\beta_j\cp_j$, so $s_{ij}=\beta_jv_{ij}$.
\begin{itemize}[leftmargin=1.5em]
    \item Maximization over $p_{ij}\in[0,1]$ of $p_{ij}(s_{ij}-c_i)$ gives condition~(i).
    \item Complementary slackness $c_i(K_i-\sum_jp_{ij})=0$ with primal feasibility gives condition~(ii).
    \item Maximization over $u_j$ gives $u_j=B_j/\beta_j$; since $\beta_j>0$, complementary slackness makes the utility constraint tight, $\sum_iv_{ij}p_{ij}+\omega_j=B_j/\beta_j$. Multiplying by $\beta_j$ gives $S_j=B_j-\beta_j\omega_j\le B_j$, which is condition~(iii).
    \item Maximization over $\omega_j\ge 0$ of $(\beta_j-1)\omega_j$ gives $\omega_j=0$ whenever $\beta_j<1$. So $S_j<B_j$ implies $\omega_j>0$, hence $\beta_j=1$ and $b_j=\cp_j$: condition~(iv).
\end{itemize}

\paragraph{Equilibrium implies optimality.} Let $(\mathbf b,\mathbf p)$ be a multi-slot FPPE with thresholds $\mathbf c$. Set $\beta_j=b_j/\cp_j\in(0,1]$, $u_j=B_j/\beta_j$, and $\omega_j=(B_j-S_j)/\beta_j\ge 0$ (by (iii)). Then $\sum_iv_{ij}p_{ij}+\omega_j=S_j/\beta_j+\omega_j=u_j$, so $(\mathbf p,\mathbf u,\bm\omega)$ is feasible with tight utility constraints, and $(\bm\beta,\mathbf c)$ is dual feasible. Condition~(i) says $\mathbf p$ maximizes $\mathcal G$ in $\mathbf p$; $u_j=B_j/\beta_j$ maximizes it in $\mathbf u$; and by (iv), $\beta_j<1$ implies $S_j=B_j$, i.e., $\omega_j=0$, so $\bm\omega$ maximizes it in $\bm\omega$. Hence $\mathcal G(\mathbf p,\mathbf u,\bm\omega;\bm\beta,\mathbf c)$ equals the dual objective at $(\bm\beta,\mathbf c)$. By tight utility constraints and condition~(ii), the multiplier terms of $\mathcal G$ vanish, so it also equals the primal objective. Equal primal and dual objectives certify optimality of both by weak duality.

\paragraph{Consequences.} Existence follows because optimal primal and dual solutions exist. By Theorem~\ref{thm:hard-dual}, $\bm\beta$ is dual-optimal if and only if $\mathbf x=\log(\bm\beta\odot\bm{\cp})$ minimizes $\Lhard$ over $\Xinf$, which identifies FPPE bids with minimizers of $\Lhard$. \hfill$\square$

\subsection{Proof of Theorem~\ref{thm:unique}}\label{app:unique}
Work in bid coordinates. Let $G(\mathbf b):=\Lhard(\log\mathbf b)=\sum_i\TopK_{K_i}(\mathbf b\odot\mathbf r_i)-\sum_jB_j\log b_j$ on the box $\{0<b_j\le\cp_j\}$. Each $\mathbf b\mapsto\TopK_{K_i}(\mathbf b\odot\mathbf r_i)$ is a convex function composed with a linear map, hence convex. Each $b_j\mapsto -B_j\log b_j$ is strictly convex, since its second derivative is $B_j/b_j^2>0$. Hence $G$ is strictly convex on a convex set and has at most one minimizer; one exists by Theorem~\ref{thm:hard-dual}. As $\mathbf x\mapsto e^{\mathbf x}$ is a bijection from $\Xinf$ onto the box, $\mathbf x^\star=\log\mathbf b^\star$ is the unique minimizer of $\Lhard$. The last claim follows from Theorem~\ref{thm:characterization}. \hfill$\square$

\subsection{Proof of Corollary~\ref{cor:essential}}\label{app:essential}
(a) is Theorem~\ref{thm:unique} and (b) is Definition~\ref{def:fppe}(iv). For (c) and (d), take the thresholds to be optimal dual prices (Theorem~\ref{thm:characterization}); for fixed $\bm\beta^\star$ they minimize the per-request term of \eqref{eq:dk}, so $c_i\in[s^\star_{i,(K_i+1)},s^\star_{i,(K_i)}]$ by Lemma~\ref{lem:topk}. If $c_i=0$, all scores (which are positive) exceed $c_i$, so every candidate wins and $K_i=M_i$; revenue is $\sum_js^\star_{ij}=\TopK_{K_i}(\mathbf s^\star_i)$. If $c_i>0$, capacity binds: campaigns scoring above $c_i$ win with probability one, campaigns scoring below lose, and the remaining capacity is shared by campaigns scoring exactly $c_i$, each paying $c_i$. Because $c_i\in[s^\star_{i,(K_i+1)},s^\star_{i,(K_i)}]$, at most $K_i$ scores exceed $c_i$ and at least $K_i$ scores are $\ge c_i$, so the purchased mass is exactly the multiset of the top $K_i$ scores and revenue is $\TopK_{K_i}(\mathbf s^\star_i)$.

For (d), suppose $s^\star_{i,(K_i)}>s^\star_{i,(K_i+1)}$. If $c_i<s^\star_{i,(K_i)}$, every top-$K_i$ campaign scores strictly above $c_i$ and wins with probability one, and capacity forces every other campaign to zero. If $c_i=s^\star_{i,(K_i)}$, then $c_i>0$, every other campaign scores strictly below $c_i$ and loses, and binding capacity forces all $K_i$ top campaigns to win with probability one.

\section{Deferred Proofs for Section~\ref{sec:fd_relaxation}}\label{app:proofs-smooth}

\subsection{Proofs of Lemma~\ref{lem:fd-dual} and Proposition~\ref{prop:reg-market}}\label{app:fd}

\paragraph{Per-coordinate identity.} For $z\in\R$ and $T>0$, the function $p\mapsto pz+T\Hb(p)$ is strictly concave on $[0,1]$ with derivative $z-T\log\frac{p}{1-p}$, which vanishes at $p=\sigma(z/T)\in(0,1)$. Substituting, $\sigma(z/T)z+T\Hb(\sigma(z/T))=T\log(1+e^{z/T})$.

\paragraph{Proof of Lemma~\ref{lem:fd-dual}.} Dualize the single equality constraint $\sum_jp_{ij}=K_i$ of \eqref{eq:primal} with a free multiplier $c$. For fixed $c$, the Lagrangian separates across $j$, and the per-coordinate identity gives the dual function $K_ic+T_i\sum_j\log(1+e^{(s_{ij}-c)/T_i})$. The problem maximizes a concave function subject to one affine equality over $[0,1]^{M_i}$, and for $K_i<M_i$ the point $p_{ij}=K_i/M_i$ is feasible and lies in the relative interior of the box, so strong duality holds~\citep[Section~5.2.3]{boyd2004convex}, giving \eqref{eq:Vi}. The dual objective is strictly convex in $c$ with derivative $K_i-F_i(c)$, so its unique minimizer is the root of Theorem~\ref{thm:sigmoid}. The primal maximizer is unique by strict concavity of $\Hb$, and maximizing the Lagrangian at the optimal $c_i$ gives \eqref{eq:sigmoid}. If $K_i=M_i$, the only feasible point is $\mathbf p_i=\mathbf 1$, whose entropy is zero.

\paragraph{Proof of Proposition~\ref{prop:reg-market}.} Requests with $K_i=M_i$ have $p_{ij}\equiv 1$ and contribute constants; remove them. Dualize the utility constraints with $\beta_j\ge 0$ and the capacity equalities with free $c_i$. The suprema over $u_j$ and $\omega_j$ are as in Section~\ref{ssec:lagrangian}; the supremum over $\mathbf p_i\in[0,1]^{M_i}$ of $\sum_j[p_{ij}(\beta_jv_{ij}-c_i)+T_i\Hb(p_{ij})]$ is $T_i\sum_j\log(1+e^{(\beta_jv_{ij}-c_i)/T_i})$ by the per-coordinate identity. Minimizing over $c_i$ gives $V_i(\bm\beta\odot\mathbf v_i;T_i)$ by Lemma~\ref{lem:fd-dual}, and the change of variables of Theorem~\ref{thm:hard-dual} gives $\Lsm(\mathbf x)+C_0$. The point $p_{ij}=K_i/M_i$, $\omega_j$ large, $u_j$ small is feasible and in the relative interior of the domain, so the refined Slater condition gives zero duality gap and dual attainment; primal attainment follows as in Appendix~\ref{app:slater}, the regularizer being bounded. The saddle-point conditions give the sigmoid allocation, $u_j=B_j/\beta_j$ and $(1-\beta_j)\omega_j=0$, from which $S_j=B_j-\beta_j\omega_j$ yields conditions (iii)--(iv) exactly as in Appendix~\ref{app:characterization}. \hfill$\square$

\subsection{Proof of Lemma~\ref{lem:gap}}\label{app:gap}
The lower bound holds because the hard top-$K_i$ indicator lies in $\mathcal C_{K_i}$ and has zero entropy. For the upper bound, $\sum_jp_{ij}s_{ij}\le\TopK_{K_i}(\mathbf s_i)$ on $\mathcal C_{K_i}$, and by concavity of $\Hb$ (Jensen's inequality) $\sum_j\Hb(p_{ij})\le M_i\Hb(\frac{1}{M_i}\sum_jp_{ij})=M_i\Hb(K_i/M_i)$. Finally, $\Hb\le\log 2$, and $M\Hb(K/M)=K\log\frac{M}{K}+(M-K)\log\bigl(1+\frac{K}{M-K}\bigr)\le K\log\frac MK+K$.

\subsection{Proofs of Theorem~\ref{thm:sigmoid} and Corollary~\ref{cor:bracket}}\label{app:sigmoid}
\paragraph{Theorem~\ref{thm:sigmoid}.} Each summand is continuous and strictly decreasing in $c$ with values in $(0,1)$, tending to $1$ as $c\to-\infty$ and to $0$ as $c\to+\infty$; this gives monotonicity, the derivative, and the range. Since $0<K_i<M_i$, the intermediate value theorem gives a root, unique by strict monotonicity; the implicit function theorem gives smoothness because $F_i'\neq 0$. The derivative of the objective in \eqref{eq:Vi} with respect to $c$ is $K_i-F_i(c)$, so $c_i^\star$ is its minimizer. For the bracket, note that $\sigma(-\rho_i/T_i)=K_i/M_i$. Since $s_{ij}\ge\min_ks_{ik}$ for all $j$, $F_i(\min_ks_{ik}+\rho_i)\ge M_i\,\sigma(-\rho_i/T_i)=K_i$; since $s_{ij}\le\max_ks_{ik}$, $F_i(\max_ks_{ik}+\rho_i)\le K_i$. Monotonicity places the root in between.

\paragraph{Corollary~\ref{cor:bracket}.} Sufficiency follows from \eqref{eq:bracket}. For necessity, if all scores equal $s$ then $F_i(c)=M_i\sigma((s-c)/T_i)$ and the root is exactly $s+\rho_i$, which lies in $[s-\kappa T_i,s+\kappa T_i]$ only if $|\rho_i|\le\kappa T_i$.

\subsection{Proof of Lemma~\ref{lem:mismass}}\label{app:mismass}
$U_i=D_i$ because $\sum_jp_{ij}=K_i$. Let $m$ be the midpoint of the gap. If $c_i\ge m$, each candidate outside the top $K_i$ has $s_{ij}-c_i\le-g_i/2$, so $U_i\le(M_i-K_i)\sigma(-g_i/(2T_i))$. If $c_i\le m$, each top-$K_i$ candidate has $s_{ij}-c_i\ge g_i/2$, so $D_i\le K_i\sigma(-g_i/(2T_i))$. Finally $\sigma(-z)\le e^{-z}$.

\subsection{Proof of Theorem~\ref{thm:gradient}}\label{app:gradient}
For $K_i<M_i$, $V_i$ is the value of a minimization over $c$ whose objective is smooth in $(c,\mathbf s_i)$ and has the unique minimizer $c_i^\star(\mathbf s_i)$, itself smooth (Theorem~\ref{thm:sigmoid}). By the envelope theorem (Danskin's theorem), $\partial V_i/\partial s_{ij}$ equals the partial derivative of the objective at $c_i^\star$, namely $\sigma((s_{ij}-c_i^\star)/T_i)=p_{ij}$; the chain rule gives \eqref{eq:grad}. For $K_i=M_i$, $V_i=\sum_js_{ij}$ and $p_{ij}=1$. Summing over $i$ and adding the derivative of $-B_jx_j$ gives \eqref{eq:SminusB}.

\subsection{Proofs of Propositions~\ref{prop:convex_smooth} and~\ref{prop:elasticity}}\label{app:convex}
\textbf{(ii) Smoothness.} The root $c_i^\star(\mathbf s_i)$ is $C^\infty$ by Theorem~\ref{thm:sigmoid}, hence so is $V_i(\mathbf s_i;T_i)=c_i^\star K_i+T_i\sum_j\log\bigl(1+e^{(s_{ij}-c_i^\star)/T_i}\bigr)$; for $K_i=M_i$, $V_i$ is linear. Composing with $\mathbf s_i(\mathbf x)=\mathbf r_i\odot e^{\mathbf x}$ and adding the affine term, $\Lsm\in C^\infty$.

\textbf{(i) Convexity.} For each $i$, $V_i(\mathbf s;T_i)$ is the pointwise supremum of a family of affine functions of $\mathbf s$, hence convex, and $\partial V_i/\partial s_j=p_{ij}\in[0,1]$, so $V_i$ is non-decreasing in each coordinate. The map $x_j\mapsto r_{ij}e^{x_j}$ is convex, so $\mathbf x\mapsto V_i(\mathbf s_i(\mathbf x);T_i)$ is convex, and so is $\Lsm$. For the sandwich, differentiate \eqref{eq:SminusB} using \eqref{eq:dpds}: for $\mathbf y\in\R^M$,
\begin{equation*}
    \mathbf y^\top\nabla^2\Lsm\,\mathbf y=\sum_jS_jy_j^2+\sum_i\frac{1}{T_i}\Biggl[\sum_ja_{ij}s_{ij}^2y_j^2-\frac{\bigl(\sum_ja_{ij}s_{ij}y_j\bigr)^2}{\sum_ja_{ij}}\Biggr],
\end{equation*}
where requests with $K_i=M_i$ contribute nothing to the second sum. By the Cauchy--Schwarz inequality the bracket lies in $\bigl[0,\sum_ja_{ij}s_{ij}^2y_j^2\bigr]$, and $a_{ij}\le\frac14$. This gives \eqref{eq:hessian-sandwich}; $S_j>0$ because every campaign is eligible somewhere and $p_{ij}>0$.

\textbf{(iii) Gradient Lipschitz on compacta.} On a compact convex $\mathcal Q$ the scores and spends are bounded, and the upper bound in \eqref{eq:hessian-sandwich} bounds the operator norm of the Hessian by the stated constant; since $\mathcal Q$ contains the segment between any two of its points, the mean value theorem turns this into the Lipschitz bound. Each $1/T_i$ multiplies bounded factors (the scores and the Bernoulli variances), so $L_\mathbf T=O(\max_i1/T_i)$.

\textbf{Proposition~\ref{prop:elasticity}.} Differentiate $S_j=\sum_ip_{ij}s_{ij}$ using $\partial s_{ij}/\partial x_k=\delta_{jk}s_{ij}$ and \eqref{eq:dpds}; then use $w_{ij}\le 1$.

\subsection{Proof of Lemma~\ref{lem:floor}}\label{app:floor}
If $b_j<\cp_j$ at a minimizer, then $S_j=B_j$ (Definition~\ref{def:fppe}(iv) in the hard case, Proposition~\ref{prop:reg-market} in the smoothed case), while $S_j=b_j\sum_ip_{ij}r_{ij}\le b_jR_j$; so $b_j\ge B_j/R_j$. The minimizer over $\Xinf$ therefore lies in $\mathcal X$ and minimizes over $\mathcal X$ as well; it is the only minimizer over $\mathcal X$ by the strict convexity used in Theorem~\ref{thm:unique} and Proposition~\ref{prop:smooth-unique}.

\subsection{Proof of Proposition~\ref{prop:smooth-unique}}\label{app:smooth-unique}
Uniqueness: $\Lsm$ is strictly convex (Proposition~\ref{prop:convex_smooth}(i)). Existence over $\Xinf$: by Lemma~\ref{lem:gap}, $V_i\ge\TopK_{K_i}\ge 0$, and $x_k\le\log\cp_k$ on $\Xinf$, so for every $j$,
\begin{equation*}
    \Lsm(\mathbf x)\ \ge\ -\sum_kB_kx_k\ \ge\ -B_jx_j-\sum_{k\ne j}B_k\log\cp_k.
\end{equation*}
Every coordinate is therefore bounded below on each sublevel set, and bounded above by $\log\cp_j$; the sublevel sets are closed and bounded, and the continuous $\Lsm$ attains its minimum. Existence over the compact box $\mathcal X$ is immediate.

\subsection{Proof of Theorem~\ref{thm:convergence}}\label{app:convergence}
\textbf{Step 1 (one-sided approximation).} By Lemma~\ref{lem:gap}, uniformly in $\mathbf x$,
\begin{equation}\label{eq:onesided}
    0\ \le\ \Lsm(\mathbf x)-\Lhard(\mathbf x)\ \le\ \Gamma_\mathbf T.
\end{equation}
\textbf{Step 2 (near-optimality).} By \eqref{eq:onesided} and minimality of $\mathbf x^\star(\mathbf T)$,
\begin{equation*}
    \Lhard\bigl(\mathbf x^\star(\mathbf T)\bigr)\le\Lsm\bigl(\mathbf x^\star(\mathbf T)\bigr)\le\Lsm(\mathbf x^\star)\le\Lhard(\mathbf x^\star)+\Gamma_\mathbf T,
\end{equation*}
so $\mathbf x^\star(\mathbf T)$ is a $\Gamma_\mathbf T$-minimizer of the hard potential.

\textbf{Step 3 (strong convexity in bid coordinates).} Let $G(\mathbf b)=\Lhard(\log\mathbf b)$ as in the proof of Theorem~\ref{thm:unique}, and let $\|\mathbf y\|_W^2:=\sum_j(B_j/\cp_j^2)y_j^2$. On $(0,\cp_j]$ the second derivative of $-B_j\log b_j$ is $B_j/b_j^2\ge B_j/\cp_j^2$, so $G-\frac12\|\cdot\|_W^2$ is convex on the box. For $\mathbf b$ in the box and the midpoint $\mathbf m=(\mathbf b+\mathbf b^\star)/2$,
\begin{equation*}
    G(\mathbf b^\star)\le G(\mathbf m)\le\tfrac12G(\mathbf b)+\tfrac12G(\mathbf b^\star)-\tfrac18\|\mathbf b-\mathbf b^\star\|_W^2,
\end{equation*}
so $\|\mathbf b-\mathbf b^\star\|_W^2\le 4\bigl(G(\mathbf b)-G(\mathbf b^\star)\bigr)$. Taking $\mathbf b=\mathbf b^\star(\mathbf T)$ and using Step~2 gives \eqref{eq:rate}; the bound on $\Gamma_\mathbf T$ is Lemma~\ref{lem:gap}.

\paragraph{What the bound does and does not say.} Dividing \eqref{eq:rate} by $\sum_jB_j$, the budget-weighted root-mean-square error in pacing multipliers is at most $2\sqrt{\Gamma_\mathbf T/\sum_jB_j}$. The bound is informative when the total smoothing gap is small relative to total budget, that is, when each $T_i$ is small relative to the scores it smooths; this is what the schedule $T_i=\alpha_{\text{temp}}\cdot\mathrm{std}(\mathbf s_i)$ with $\alpha_{\text{temp}}\to 0.01$ aims for (Section~\ref{sec:bid_update}). The bound is worst-case, and we do not claim it is tight.

\paragraph{Allocations.} Combining Theorem~\ref{thm:convergence} with Lemma~\ref{lem:mismass}: at every request whose FPPE boundary gap $g_i^\star=s^\star_{i,(K_i)}-s^\star_{i,(K_i+1)}$ is positive, the scores at $\mathbf b^\star(\mathbf T)$ have gap at least $g_i^\star/2$ once $\mathbf T$ is small, so the soft allocation converges to the unique hard allocation of Corollary~\ref{cor:essential}(d), with mis-assigned mass at most $M_ie^{-g_i^\star/(4T_i)}$. Mass stays soft only at requests that tie, or nearly tie, at the boundary.

\paragraph{Why the argument is short.} The key structural fact is one-sided: the softplus upper-bounds the hinge, so $\Lsm\ge\Lhard$ everywhere. This gives near-optimality of $\mathbf x^\star(\mathbf T)$ without a specific annealing schedule; uniqueness of the hard minimizer prevents different annealing paths from converging to different limits; and strong convexity of the logarithmic regularizer in bid coordinates turns near-optimality into a rate. Together these avoid the general machinery of epi-convergence~\citep{rockafellar1998variational}. By Lemma~\ref{lem:floor}, everything above holds verbatim on the box $\mathcal X$ once $\varepsilon\le\min_j\underline b_j$.

\section{Machine-Checked Proofs}\label{app:lean}

Every lemma, proposition, theorem, and corollary of Sections~\ref{sec:model}--\ref{sec:practice}, together with the definitions it uses, is formalized and machine-checked in Lean~4~\citep{demoura2021lean4} with the Mathlib library~\citep{mathlib2020} and the optimization library Optlib~\citep{li2025optimality}. So are three claims made in passing: the first-order conditions on the box and their sufficiency (Section~\ref{sec:kkt}), the paragraph on allocations in Appendix~\ref{app:convergence}, and item (S1) of Section~\ref{sec:scope}, for any update direction that is continuous in the bids. The Lean kernel checks every proof; the formalization contains no unproved steps and uses only Lean's three standard axioms. Two things are not formalized: the convergence of projected gradient descent invoked for Algorithm~\ref{alg:solver}, which we cite rather than prove, and the asymptotic remarks on the Lipschitz rate.

\paragraph{Conventions.} The potentials are stated in bids $b_j=e^{x_j}$ on the box $0<b_j\le\cp_j$, the image of $\Xinf$ under $\mathbf x\mapsto e^{\mathbf x}$. The formal model lets a campaign appear more than once among a request's candidates, so unit demand is a special case; Proposition~\ref{prop:elasticity} is proved in the corresponding general form, which reduces to the stated formulas under unit demand. The bounds \eqref{eq:hessian-sandwich} are proved as bounds on $\mathbf y^\top\nabla^2\Lsm(\mathbf x)\,\mathbf y$ for every $\mathbf y$, and strict convexity, the bound $\eta_j\ge 1$, item (S1), and the existence of an admissible floor in Lemma~\ref{lem:floor} use the standing assumption that every campaign is eligible for some request.

\paragraph{Proofs that differ.} Strong duality (Proposition~\ref{prop:slater}) is obtained from primal attainment by compactness and the Karush--Kuhn--Tucker conditions under a linear constraint qualification, which Optlib provides, rather than from Slater's condition. The gradient identity (Theorem~\ref{thm:gradient}) is proved without the envelope theorem, by squeezing $V_i$ between its value at a fixed allocation and the threshold objective at the fixed root. Smoothness of the root (Theorem~\ref{thm:sigmoid}) follows from the inverse function theorem applied to $(\mathbf s,T,c)\mapsto(\mathbf s,T,F_i(c))$.

\section{Why Second-Price Breaks the Potential Structure}\label{app:spa}

This appendix expands on the first-price argument of Section~\ref{sec:intro}; SPA and FPA denote second- and first-price auctions. Let $S_j(b_j)$ be campaign $j$'s expected spend as a function of its own bid, all other bids fixed. In FPA a winner pays its own bid, so $S_j(b_j)=\sum_i p_{ij}(b_j)\,s_{ij}(b_j)$, where $p_{ij}$ is the probability that $j$ wins a slot of request $i$ and $s_{ij}=b_j r_{ij}$ its expected payment when it does; both factors are non-decreasing in $b_j$. Hence $S_j(b_j)$ is (1) \textbf{non-decreasing} in $b_j$, with elasticity $\partial\log S_j/\partial\log b_j$ at least one under the entropic relaxation (Proposition~\ref{prop:elasticity}); (2) unbounded as $b_j\to\infty$, under mild coverage of the request pool; and (3) \textbf{smooth} under the relaxation. These per-campaign properties reflect a global structure: the equilibrium minimizes a convex potential whose gradient is the spend--budget mismatch (Theorem~\ref{thm:gradient}), the minimizer is unique for the hard market (Theorem~\ref{thm:unique}) and every relaxation (Proposition~\ref{prop:smooth-unique}), and the relaxed potential is smooth with curvature $O(1/T)$ (Proposition~\ref{prop:convex_smooth}). This is what lets our gradient-based solver reach the equilibrium of the relaxed problem rather than merely approximate it.

SPA breaks all three properties for one reason: a winner pays the runner-up's score, which is \emph{decoupled} from its own bid. Below, (a) the gradient identity has no SPA analogue we know of, (b) spend is flat in a winner's own bid over wide ranges and jumps when bids cross, so the signal budget pacing relies on is mostly absent, and (c) smoothing the order statistic that sets SPA prices is expected to break convexity. The literature documents the same contrast for equilibria: second-price pacing equilibria always exist but can be multiple~\citep{conitzer2022multiplicative} and are PPAD-hard to compute~\citep{chen2024complexity}, whereas the first-price equilibrium is unique and solves a convex program~\citep{conitzer2022pacing}; see also~\citep{edelman2007gsp,balseiro2019learning}. This is why we model the market as first-price: the choice concerns equilibrium computation only, not revenue, incentives, or the mechanism a platform deploys.

\paragraph{(a) The gradient identity breaks.} In FPA the price on a won request is the own score $s_{ij}$, so differentiating $V_i$ with respect to $x_j=\log b_j$ gives exactly $p_{ij}\,s_{ij}$, the expected spend contribution, and summing over $i$ gives $\nabla_{x_j}\Lsm=S_j-B_j$ \eqref{eq:SminusB}: a convex objective with exact spend-mismatch gradient, the identity that drives this paper. In SPA the per-win payment is the runner-up's score, which depends on another campaign's $x_k$, so the direct dependence on $x_j$ no longer evaluates to the expected spend and the envelope argument does not apply. We know of no convex potential whose gradient is the SPA spend mismatch, though we make no formal impossibility claim.

\paragraph{(b) Zero-gradient plateaus and steps.} Suppose campaign $j$ wins a request at score $s_j=10$ with runner-up $s_k=5$. Raising its score to 12 or lowering it to 8 leaves the payment at 5: $\partial S_j/\partial b_j=0$ wherever $j$ stays the winner, and spend changes only when $b_j$ crosses a neighbor's bid, a measure-zero event in the hard limit. A continuous relaxation turns this into nearly flat plateaus joined by sharp transitions, so the signal that budget pacing relies on is mostly absent.

\paragraph{(c) Smoothing order statistics is expected to break convexity.} A consistent SPA relaxation must smooth not the top-$K$ selection, as our threshold sigmoid does for FPA, but the order statistic that sets the price: for a single slot, the second-largest score, the highest among a winner's competitors. Differentiable order statistics exist (NeuralSort~\citep{grover2019neuralsort} and Sinkhorn-based operators~\citep{cuturi2019differentiable,xie2020differentiable}), but the resulting expected-spend surrogate is generally not expected to be convex: lowering a bid can leave spend unchanged over a range and then drop it abruptly when the winner loses the slot. We do not prove non-convexity, but the plateau-and-step geometry of (b) should persist qualitatively through standard order-statistic smoothers, and a quasi-Newton method such as L-BFGS-B~\citep{byrd1995lbfgsb} would then meet negative-curvature secant pairs or uninformative gradients on the plateaus, corrupting its inverse-Hessian approximation.

\section{Implementation Details}\label{app:impl}

\subsection{Hyperparameters}\label{app:hyper}
Table~\ref{tab:hyper} lists the constants of the main configuration (Section~\ref{sec:exp-setup}).

\begin{table}[htbp]
\centering
\small
\caption{Hyperparameters of the main configuration.}
\label{tab:hyper}
\begin{tabular}{llr}
\toprule
Component & Symbol & Value \\
\midrule
Temperature coefficient: start, end & $\alpha_0$, $\alpha_{\text{end}}$ & 1.25, 0.01 \\
Anneal length (geometric decay) & $t_{\text{horizon}}$ & 80 iterations \\
Step coupling, $\alpha_{\text{step}}=c\,\alpha_{\text{temp}}$ & $c$ & 0.42 \\
Ratio clamp & $[r_{\mathrm{lo}},r_{\mathrm{hi}}]$ & $[0.1,10]$ \\
Bid floor & $\varepsilon$ & $10^{-8}$ \\
Temperature floor & $\min_iT_i$ & $10^{-6}$ \\
Momentum coefficients & $\gamma_t$ & FISTA sequence \\
Momentum warm-up & & 3 iterations \\
Momentum restart & & adaptive; and at end of anneal \\
Root-solve bracket margin, in units of $T_i$ & $\kappa$ & 50 \\
Bisection, then Newton iterations & & 15, then 5 \\
Logit clamp & & $\pm 50$ \\
Iteration budget & & 800 \\
\bottomrule
\end{tabular}
\end{table}

\subsection{Numerical stability of the root-finding kernel}\label{app:numerics}
\paragraph{Safe bracketing.} So that convergence does not depend on candidate-pool depth, we bracket $c_i\in[s_{\min}-\kappa T_i,\ s_{\max}+\kappa T_i]$ with $\kappa=50$, where $s_{\min},s_{\max}$ are request $i$'s smallest and largest eligible scores: every logit is $\ge 50$ at the lower bound ($\sigma\approx 1-2\cdot 10^{-22}$, $F_i\approx M_i\gg K_i$) and $\le -50$ at the upper bound ($\sigma\approx 2\cdot 10^{-22}$, $F_i\approx 0<K_i$ for any realistic candidate pool). By Corollary~\ref{cor:bracket}, a $\kappa T$ bracket is valid for every score vector if and only if $\kappa\ge|\log((M_i-K_i)/K_i)|$. For $\kappa=50$ this holds whenever $M_i/K_i<1+e^{50}$, always in practice; for $\kappa=5$ it fails as soon as $M_i/K_i>1+e^5\approx 149$, matching the empirical failure of a $5T$ bracket when $M_i\gtrsim150$ and $K_i$ is small, where the sum at $c_{\mathrm{hi}}$, about $M_i\,\sigma(-5)\approx 0.0067\,M_i$, can exceed $K_i$ and the root escapes the bracket. The tighter bracket \eqref{eq:bracket} of Theorem~\ref{thm:sigmoid} has width $s_{\max}-s_{\min}$ instead of $s_{\max}-s_{\min}+2\kappa T_i$ and would save $\log_2\bigl(1+2\kappa T_i/(s_{\max}-s_{\min})\bigr)$ bisection steps; our kernel does not yet use it.

\paragraph{Hybrid bisection--Newton solver.} We use a two-phase vectorized solver:
\begin{enumerate}
    \item \textbf{Bisection (15 iterations)}: halves the bracket, preserving $F_i(c_{\mathrm{lo}})>K_i>F_i(c_{\mathrm{hi}})$. After 15 halvings the width is $2^{-15}(s_{\max}-s_{\min}+100\,T_i)$, a relative reduction of $3\times 10^{-5}$ that stays well above \texttt{float32} epsilon ($1.2\times 10^{-7}$), leaving numerical headroom for Newton.
    \item \textbf{Newton polish (5 iterations)}: we apply Newton's method to $F_i(c)-K_i$, using $F_i'(c)=-\tfrac{1}{T_i}\sum_j p_{ij}(1-p_{ij})<0$ (Theorem~\ref{thm:sigmoid}), which gives the update
    \begin{equation}
        c^{(t+1)} \;=\; c^{(t)} \;+\; \frac{\sum_j p_{ij}-K_i}{(1/T_i)\sum_j p_{ij}(1-p_{ij})}.
    \end{equation}
    If $\sum_j p_{ij}>K_i$ the threshold is too low, and the positive numerator raises it. Convergence is quadratic in the basin. The step is clamped to the final bisection bracket, because in the saturated regime the derivative approaches zero and can amplify \texttt{float32} quantization noise in $\sum_j p_{ij}-K_i$ into an explosive step.
\end{enumerate}
In the fused kernels (Section~\ref{sec:solver}) both phases run inside each request's program; the earlier implementation vectorized them across all requests with segmented reductions over the flattened list of $(i,j)$ candidate pairs. The kernel also applies a logit clamp at $\pm50$, the exact $K_i=M_i$ short-circuit of Lemma~\ref{lem:fd-dual}, an empty-row guard, a variance clamp, and a post-solve renormalization, listed below.

\begin{itemize}
    \item \textbf{Logit clamp $\pm 50$} before the sigmoid to prevent \texttt{float32} overflow; $\sigma(50)$ rounds to $1$, and $\sigma(-50)\approx 2\times 10^{-22}$ is negligible.
    \item \textbf{$K_i=M_i$ short-circuit}: every candidate wins, but saturated sigmoids cannot represent $p=1$ exactly at large $M_i$, so we set $p\leftarrow 1$, the exact solution (Lemma~\ref{lem:fd-dual}).
    \item \textbf{Empty-row guard}: rows filtered upstream leave $s_{\min}=+\infty$, $s_{\max}=-\infty$; we set $c=0$ (empty rows contribute nothing downstream).
    \item \textbf{Variance clamp $\ge 0$} before the square root to guard against \texttt{float32} rounding producing tiny negatives.
    \item \textbf{Post-solve renormalization}: $\sum_j p_{ij}$ may differ from $K_i$ by $O(10^{-6})$ because the bracket does not collapse to machine epsilon, so we rescale $p\leftarrow \mathrm{clamp}(p\cdot K_i/\sum p,\ 0,\ 1)$. When the factor $s$ exceeds 1 the clamp can leave the sum slightly \emph{below} $K_i$ (by at most $\sim 10^{-5}$), a loss of exact capacity we accept to keep $p\in[0,1]$.
\end{itemize}

\section{Allocation Audit}\label{sec:kernel_correctness}

\begin{figure}[!htbp]
\centering
\includegraphics[width=\textwidth]{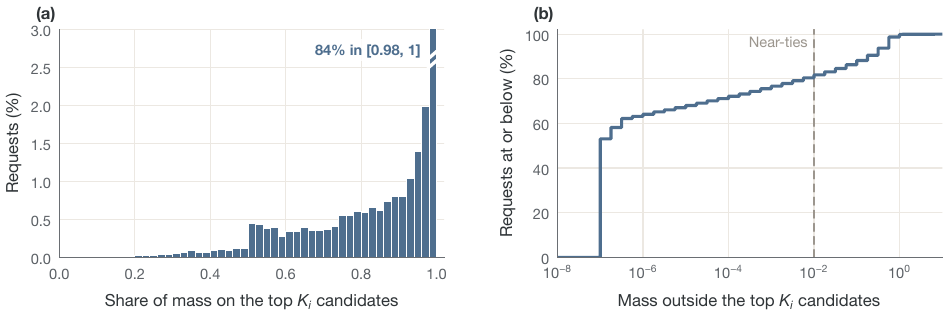}
\caption{Allocation audit on the headline day. (a) Each request's share of mass on its top $K_i$ candidates (axis cut at 3\%). (b) Cumulative share of requests by mass outside the top $K_i$; the step near $10^{-7}$ is \texttt{float32} rounding.}
\label{fig:topk}
\end{figure}

The fused kernels write no per-candidate arrays, so we audit the allocation with the earlier, unfused kernel, which computes the same smoothed allocation and can record per-candidate marginals. We rerun the main configuration on the headline day with the unfused kernel for 800 iterations and, at the last iteration, record every candidate of a uniform random 0.1\% of requests: its score, its marginal $p_{ij}$, and the request's threshold, temperature, and renormalization factor. The final bids of this run and of the fused four-GPU run differ by a median of $1.3\times10^{-6}$ and at most $9.9\times10^{-5}$ over active campaigns, and the same campaigns end at their caps; the fused runs on two and four GPUs give bit-identical final bids. On the audit sample:
\begin{itemize}
    \item \textbf{Capacity invariant}: $|\sum_j p_{ij}-K_i|$ has a maximum of $2.1\times 10^{-5}$ and a mean of $8.6\times 10^{-8}$; no request exceeds $10^{-4}$.
    \item \textbf{Top-$K$ concentration}: 64.5\% of requests place mass $K_i$ on their $K_i$ highest-scored candidates to within $10^{-6}$, and the median request is fully concentrated (Figure~\ref{fig:topk}). For 18.2\% of requests more than $0.01$ of the mass lies outside the top $K_i$; the 99th percentile of that mass is $0.60$.
    \item \textbf{Near-ties}: these soft requests are near-ties at the boundary. Their gap $g_i$ between the $K_i$-th and $(K_i{+}1)$-th scores is a median of 3.1 temperatures $T_i$ and at most 11.4, whereas the median fully concentrated request has a gap of 94 temperatures. The bound of Lemma~\ref{lem:mismass}, $\sum_{j\notin\text{top-}K}p_{ij}\le\max\{K_i,M_i-K_i\}\,e^{-g_i/(2T_i)}$, holds on every audited request. On a near-tie the hard top-$K$ solution is itself nearly non-unique, and the smoothed allocation splits the mass among the tied candidates.
    \item \textbf{Post-solve renormalization} is non-trivial ($|s-1|\ge 10^{-6}$) for 9.1\% of requests; when active, $s\in[0.99998, 1.00002]$.
    \item \textbf{Gradient sanity check} (fixed-$\mathbf{T}$ ideal): on a 3-request, 4-candidate unit-test fixture with fixed per-request $T_i$ and none of the stabilizations (no adaptive rescaling, logit clamp, $K_i=M_i$ short-circuit, or renormalization), central finite differences of $\Lsm(\mathbf{x})$ with step $10^{-4}$ match the prediction $S_j(\mathbf{x};\,\mathbf{T})-B_j$ to $\sim\!10^{-6}$, consistent with \texttt{float32} arithmetic, and the implicit-function formula \eqref{eq:dpds} for $\partial p/\partial s$ matches PyTorch autograd to $1.5\times 10^{-6}$. This confirms $\nabla_{x_j}\Lsm=S_j-B_j$ at the ideal-model level; the implementation adds the perturbations (S1)--(S2) of Section~\ref{sec:scope}.
\end{itemize}

\end{document}